\documentstyle{amsppt}
\magnification=\magstep 1
\TagsOnRight
\NoBlackBoxes

\vsize 23truecm\hsize 16.5truecm

J.  Anal. Math. (to appear)

\vbox{\vskip 3truecm}
\parindent6mm
\voffset=-0.5truecm

\topmatter
\title  Asymptotic behavior of polynomials orthonormal on a homogeneous set
\endtitle
\rightheadtext{  Asymptotic behavior of orthonormal polynomials}
\author Franz Peherstorfer and  Peter Yuditskii \endauthor

 \address  
 Institute for Analysis and Computational
 Mathematics, Johan\-nes Kepler University of Linz,
 A--4040 Linz, Austria
 \endaddress
 \email Franz.Peherstorfer\@jk.uni-linz.ac.at \endemail
 \address 
  Dept.\ of Mathematics, Michigan State University,
 East Lansing, MI 48824, USA
 \endaddress
 \email yuditski\@math.msu.edu \endemail

\abstract {Let $E$ be a homogeneous compact set, for instance  a Cantor set
of positive length.
Further let
$\sigma$ be a positive measure with $\text{supp}(\sigma)=E$. Under the
condition that the absolutely continuous part of $\sigma$ satisfies a
Szeg\"o--type condition we give an asymptotic representation, on and off the
support, for the polynomials orthonormal with respect to $\sigma$. For the
special case that $E$ consists of a finite number of intervals and that
$\sigma$ has no singular component this is  a nowaday well known result of
Widom. If $E=[a,b]$ it becomes a classical result due to Szeg\"o and in case
that there appears in addition a singular component, it is due to
Kolmogorov--Krein. In fact the results are presented for the more general case
that the orthogonality measure may have a denumerable set of mass--points
outside of $E$ which are supposed to accumulate on $E$ only and to satisfy
(together with the zeros of the associated Stieltjes function) the
free--interpolation Carleson--type condition. Up to the case of a finite number
of mass points this is even new for the single interval case. Furthermore, as a
byproduct of our representations, we obtain that the recurrence coefficients of
the orthonormal polynomials behave asymptotically almost periodic. Or in other
words the Jacobi matrices associated with the above discussed orthonormal
polynomials are compact perturbations  of a one--sided restriction of  almost
periodic Jacobi matrices with homogeneous spectrum. Our main tool is a theory
of Hardy spaces of character--automorphic functions and forms on Riemann
surfaces of Widom type, we use also some ideas of scattering theory for
one--dimensional Schr\"odinger equations.}
\endabstract

\thanks{This work was supported
by the Austrian Science Fund FWF, project--number P12985--TEC
}\endthanks

\endtopmatter

\document

\head Introduction
\endhead

Let 
$\sigma$ be a positive measure with a compact support.
 By
$P_n(x)=P_n(x,\sigma)$ we denote the polynomial of degree
$n$ orthonormal with respect to $\sigma$, i.e.:
$$
\int P_n(x) P_m(x)\,d\sigma(x)=\delta_{n,m}.\tag 0.1
$$
It's well known that $\{P_n\}$ satisfies a three--term recurrence
relation
$$
z P_n(z)=p_n P_{n-1}(z)+q_n P_{n}(z) + p_{n+1} P_{n+1}(z),
\quad n=1, 2,\dots ,\tag 0.2
$$
with initial data
$$
p_0 P_0(z)=1,\quad zP_0(z)=q_0 P_0(z)+p_1 P_1(z).
$$

One of the main problems is to find an explicit or at least an asymptotic
representation of
the orthonormal polynomials and their  recurrence coefficients.

In the case when the support of the measure is the single interval, say
$E=[-2,2]$, this problem has been solved by Szeg\"o and Bernstein [30, 5]
rather completely. As we shall see a crucial point was Szeg\"o's discovery that
the problem transformed to the unit circle by the well known conformal
mapping
$z:\Bbb D=\{\zeta:\ \vert\zeta\vert<1\}\to \bar\Bbb C\setminus E$,
$z(\zeta)=\frac 1\zeta+\zeta$
is closely connected with an extremal problem in  the Hardy space $H^2$. 

To be able to state Szeg\"o's and our results we first will need some basic
facts on Hardy spaces and functions of bounded characteristic [12]. 

Recall
that a function $f(\zeta)$ meromorphic in $\Bbb D$ is said to be of bounded
characteristic if 
$$
\sup_{0<r<1}\left\{\int_{\Bbb T}\log^+\vert
f(rt)\vert\,dm(t)\right\}<\infty,
$$
where 
$\Bbb T=\{t:\ \vert t\vert=1\}$ and $dm$ is the Lebesgue measure on $\Bbb T$.
It can be represented as a ratio of two holomorphic functions bounded in the
disk, that is,
$$
f(\zeta)=\frac{f_+(\zeta)}{f_-(\zeta)},
\quad\sup_{\zeta\in  \Bbb D}\vert{f_\pm(\zeta)}\vert\le 1.
$$
Such functions have the representation
$$
{f_\pm(\zeta)}=\prod\frac{\overline{\zeta_l^\pm}}
{\vert\zeta_l^\pm\vert}\frac{\zeta_l^\pm-\zeta}{1-\overline{\zeta_l^\pm}\zeta}
\exp\left\{i c^\pm+\int_\Bbb T\frac{\zeta+t}{\zeta-t}\,d\tau^\pm(t)\right\},
$$
where $\zeta_l^\pm\in \Bbb D$, $c^\pm\in \Bbb R$ and $ d\tau^\pm(t)$ are
positive measures on $\Bbb T$. One can decompose $ d\tau^\pm(t)$
into the absolutely continuous $d\tau_{a.c.}^\pm(t)$ and 
singular $ d\tau_{s.}^\pm(t)$ part. The factor
$$
{f_\pm^{\text{out}}(\zeta)}=\exp\left\{ic^\pm+\int_\Bbb
T\frac{\zeta+t}{\zeta-t}\,d\tau_{a.c.}^\pm(t)\right\}
$$
is called the outer part of the function ${f_\pm(\zeta)}$.
It is defined  uniquely (up to unimodular constant) via the boundary
values of the modulus of the given function,
$d\tau_{a.c.}^\pm(t)=-\log\vert f_\pm(t)\vert\,dm(t)$. The remaining part of the
function is called the inner part. It contains the Blaschke product and the
singular component. The function $f(\zeta)$ is of Smirnov class,
(or Nevanlinna class $N^+$), if the denominator $f_-$ is an outer
function.

As usual, $H^p$ denotes the Hardy space of functions $f(\zeta)$
 analytic on $\Bbb D$ with
$$
\Vert f\Vert_p=\sup_{0<r<1}\left\{\int_{\Bbb T}\vert
f(rt)\vert^p\,dm(t)\right\}^{\frac 1 p}<\infty.
$$
Note that any function from $H^p$ is a function of Smirnov class, and
that functions of Smirnov class obey the maximum principle in the following
form: if $f$ is of Smirnov class and 
$f(t)\in L^p,\ t\in\Bbb T$, then $f\in H^p$.

Now we are ready to state Szeg\"o's result in a suitable way.
Let the transformed measure
$\sigma^*(t)$, given by
$$
2\int_{E=[-2,2]} f(x)\, d\sigma(x)=\int_{\Bbb T}f(z(t))\,d\sigma^*(t),
$$
have a decomposition of the form
$$
d\sigma^*(t)=\rho(t) d m(t)+d\sigma^*_{s.}(t),
$$
where $\sigma_{s.}^*$ is a singular measure and
$\rho(t)$ satisfies the so--called  Szeg\"o condition
$$
\int_{\Bbb T}\log\rho(t) \,dm(t)>-\infty.\tag 0.3
$$
Then $\rho$ has a representation of the form
$$
\rho(t)=\vert D(t)\vert^2\quad \text{a.e. on}\ \Bbb T,
$$
where
$$
D(\zeta):=\exp\left\{\frac 1 2\int_{\Bbb T}\log\rho(t)\frac{t+\zeta}{t-\zeta}
\,dm(t)\right\}
$$
is an outer function. Since it is a characteristic property of an outer
function, that 
$$
\text{clos}_{L^2}\{D H^\infty\}=H^2,\tag 0.4
$$ 
we have
$$
\nu(\rho):=\inf\{\int_{\Bbb T}\vert
f(t)\vert^2\rho(t)
\,dm(t):f\in H^\infty,\ \ f(0)=1\}=D^2(0)
\tag 0.5
$$
with an extremal function 
$$
f(\zeta,\rho)=\frac{D(0)}{D(\zeta)}.
$$
Using the fact that the polynomials orthonormal on the unit
circle with  respect
to the weight function $\rho$ form a so--called minimizing sequence for problem
(0.5) (see e.g. [22]) Szeg\"o has shown that uniformly on compact subsets
of $\Omega=\bar \Bbb C\setminus E$ the asymptotic relation
$$
 P_n(z(\zeta),\sigma)\sim\frac{\zeta^{-n}}{ D(\zeta)}\tag 0.6
$$
holds. By (0.6) it follows that the recurrence coefficients
have the asymptotic behavior
$$
p_n\to 1\quad\text{and}\quad q_n\to 0\quad (n\to\infty).
$$
Szeg\"o also has given an asymptotic representation of the orthonormal
polynomials on $\Bbb T$ with respect to $L^2$--norm.
\medskip

Up to the next major step in a general characterization,
due to Akhiezer--Tomchuk [1, 3, 32] and Widom [34], it took almost 40 years.
For the case of finite number of intervals, say 
$E:=[b_0,a_0]\setminus\cup_{j=1}^{N}(a_j,b_j)$ Akhiezer and Tomchuk
derived so called comparative asymptotics, that is,
the weight functions are divided into classes and if the asymptotic behavior of
the orthonormal polynomials with respect to one weight function is 
known the asymptotics with respect to all other weight functions from this
class can be given. 
 At the end of the sixties Widom discovered the
important fact that in case of several arcs and curves, moreover 
of a finite number of intervals, the extremal problem (0.5)
has to be replaced by a much more sophisticated one which we are going to
discuss next.

To see better the parallels to Szeg\"o's theory  
 let us uniformize the domain $\Omega=\bar\Bbb C\setminus E$
 by the so called universal covering map $z(\zeta)$.
Recall that $z(\zeta)$ is a meromorphic function which maps $\Bbb D$
onto
$\Omega$ and which is automorphic with respect
to the associated Fuchsian group $\Gamma$, i.e.,
$z(\gamma(\zeta))=z(\zeta),\ \forall \gamma\in\Gamma$,
and any two preimages of $z_0\in\Omega$ are $\Gamma$--equivalent, i.e.,
$
z(\zeta_1)=z(\zeta_2)\ \Rightarrow\ \exists\gamma\in\Gamma:\ 
\zeta_1=\gamma(\zeta_2)
$.
We normalize $z(\zeta)$  by the conditions
 $z(0)=\infty$, $(\zeta z)(0)>0$.

Furthermore, a character of $\Gamma$ is a complex--valued function
$\alpha:\ \Gamma\to\Bbb T$, satisfying
$$
\alpha(\gamma_1\gamma_2)=\alpha(\gamma_1)\alpha(\gamma_2)\quad
(\gamma_1,\gamma_2\in\Gamma).
$$
The characters form an Abelian compact group denoted  by $\Gamma^*$.

For a given  character $\alpha\in\Gamma^*$ define the space of
character--automorphic functions
$$
H^\infty(\Gamma,\alpha)= \{f\in H^\infty:
f(\gamma(\zeta))=\alpha(\gamma)f(\zeta),\ \forall \gamma\in\Gamma\}.
$$

The Blaschke product
$$
b(\zeta)= b(\zeta,0)=\prod_{\gamma\in\Gamma} 
 \frac{\gamma(0)-\zeta}{1-\overline{\gamma(0)}\zeta}
 \frac{\vert\gamma(0)\vert}{\gamma(0)}
$$
is called the {\it Green's function} of $\Gamma$ with respect to
the origin.
Since it is a character--automorphic function, there exists
a
$\mu\in
\Gamma^*$ such that 
$$
b(\gamma(\zeta))=
\mu(\gamma)b(\zeta).\tag 0.7
$$
Note, if $G(z)=G(z,\infty)$ denotes the Green's function of the domain
$\Omega$, then
$$
G(z(\zeta),z(0))=-\log\vert b(\zeta,0)\vert.
$$
Without loss of generality we  assume in this paper that $(bz)(0)=1$,
i.e., the capacity of $E$ is equal to 1.

\medskip

Next let us assume for a moment that $\sigma(x)$ is absolutely continuous
on $E$, i.e.,
$d\sigma(x)=\sigma'_{a.c}(x)\,dx$. Transforming this measure by the universal
$z(t)$, we get
$$\aligned
2\int_E f(x)\,\sigma'_{a.c.}(x)\,dx=&
\int_{\Bbb E}f(z(t)) \sigma'_{a.c.}(z(t))\vert z'(t)\vert 2\pi\,dm(t)\\
=&
\int_{\Bbb E}f(z(t))\rho(t)\,dm(t),
\endaligned\tag 0.8
$$
where $\Bbb E$ is a fundamental set for the action of $\Gamma$ on $\Bbb T$.
Under the conditions that $E$ consists of a finite number of intervals and that
$\rho(t)$ satisfies the generalized Szeg\"o condition (0.3)
Widom has shown that the minimum problem, $\rho\in L^1_{dm\vert\Bbb E}$,
$$
\nu(\rho,\alpha):=\inf\{\int_{\Bbb E}\vert f(t)\vert^2\rho(t)\,dm(t):\
f\in H^\infty(\Gamma,\alpha)\ \text{and}\
f(0)=1\},\tag 0.9
$$
 and its unique extremal
function $f(\zeta,\rho,\alpha)$ (which does not belong to
$H^\infty$ in general) play a crucial role in the asymptotics of polynomials
$\{P_n\}$ orthonormal on a finite number of intervals $E$. Indeed, he proved
that the minimum deviation has the following asymptotic behavior
$$
p_0^2\dots p_n^2=\min_{a_i}\int_E\vert x^n+a_1x^{n-1}+\dots+a_n\vert^2
\,d\sigma(x)\sim\nu(\rho,\mu^{n})\tag 0.10
$$
and uniformly on compact subsets of $\Omega$ the orthonormal polynomial is given
asymptotically by
$$
(p_0\dots p_n) P_n(z(\zeta),\sigma)
b^n(\zeta)
\sim
f(\zeta,\rho,\mu^{n}).\tag 0.11
$$
He gave also asymptotics on the support. Widom used spaces of multivalued
functions on
$\Omega$, but we  presented his results in a  way suitable in what follows.

\bigskip

Thus after Akhiezer--Tomchuk's, Widom's and subsequent investigations 
[4, 13, 18, 24--27, 33] the
following  natural questions arise:

a)  Does there
still hold an asymptotic representation like (0.6) respectively (0.11)
if $\bar \Bbb C\setminus E$ is infinite connected, in particular, if
$E$ is a Cantor type set?

b) Is there an analog of the Krein--Kolmogorov-Szeg\"o Theorem, that is,
does the asymptotic representation still hold if the measure has an
arbitrary singular part on $E$?

c)  How does the asymptotic representation change if outside of $E$ a
denumerable set of mass--points is added?

d) Finally, what could be said on the asymptotic behavior of the recurrence
coefficient? 

\medskip

Let us point out that so far the answer to the questions b)-- d) have been
unknown partly even for the case of a finite number of intervals.
More precisely the answer to question b) is known for the case of one
interval only [22, Sect. 2]. Concerning question c) only the case of a finite
number of mass points could be handled (see [21]).

What
concerns question d) partial answers are known (see  [4, 18, 24, 25]) when
the support of the measure
consists of a finite number of intervals. Indeed, it's well known
nowadays
that the recurrence coefficients behave asymptotically periodic, if the
harmonic measure of each component of $E$ is rational, or in other
words if the set $E$ can be described as the inverse image of a polynomial map.
If the harmonic measure of at least one of the intervals of $E$ is not rational
then it is conjectured that the recurrence coefficients behave uniformly almost
periodic in the limit. 
Recall, that a sequence of real numbers
$\{p_n\}\in l^\infty(\Bbb Z)$ is called uniformly almost
periodic if the set of sequences $\left\{\{p_{n+n_k}\},\
n_k\in \Bbb Z\right\}$ is  precompact in $l^\infty(\Bbb Z)$. The general way to
produce a sequence of this type looks as follows: let $\Cal G$ be a compact
Abelian group, and let
$f(g)$ be a continuous function on $\Cal G$, then
$$
p_n:=f(g_0+n g_1),\quad g_0,g_1\in \Cal G,\tag 0.12
$$
is an almost periodic sequence.
If $E$ consists of two intervals for rational weights with square root
singularities at the boundary points, so called Bernstein--Szeg\"o type
weights, the almost periodic recurrence coefficients can be given even
explicitly in terms     of
elliptic functions [2, 24]. In the case of several intervals for this class
of weight functions it would be still possible to give a representation in
terms of Theta--functions of several variables, using results and methods
given in [11, 31]. 

If the weight function satisfies a
Szeg\"o--type condition only we could derive from Widom's result (we will not
carry out this because it's included in our more general results below), that
the recurrence coefficients have the following representation
$$
p_n=p_n^{(1)}+p_n^{(2)},\quad
q_n=q_n^{(1)}+q_n^{(2)},\quad n\in \Bbb Z_+,
$$
where $\{p_n^{(1)}\}$, $\{q_n^{(1)}\}$ are  half--line restrictions  of
almost periodic sequences, and
$$
p_n^{(2)}\to 0,\quad q_n^{(2)}\to 0,\quad n\to\infty.\tag 0.13
$$

If we denote by $J^{(1)}$ the Jacobi matrix associated with $\{p^{(1)}_n\}$ and
$\{q^{(1)}_n\}$, i.e.,
$$
(J^{(1)} u)_n=p^{(1)}_n u_{n-1}+q^{(1)}_n u_n+p^{(1)}_{n+1}u_{n+1},\quad u\in
l^2(\Bbb Z),\tag 0.14
$$
(0.13)  says that  $J_+=J_+^{(1)}+J_+^{(2)}$
is a compact perturbation of the compression
on $l^2(\Bbb Z_+)$ of an almost periodic Jacobi matrix $J^{(1)}$.
As usual, a Jacobi matrix
is called almost periodic if the coefficient sequences are almost periodic.

In the last years
the spectral theory of almost periodic Jacobi matrices has
been studied extensively [7, 9, 23] 
in particular in connection
with integrable systems [11, 31]. 
 For  other new 
interesting developments concerning asymptotics of orthogonal polynomials
see [10, 17]. Recently a  complete description of almost periodic Jacobi
matrices with homogeneous absolutely continuous spectrum has been given in
[29] (see Theorem [29] below).

Following Carleson, we say that a compact set
$E$ is homogeneous if
there is an $\eta>0$ such that
$$
\vert(x-\delta,x+\delta)\cap E\vert\ge \eta\delta\quad \text{for all}\
0<\delta<\text {diam}\,E \quad\text {and all}\ x\in E,
$$
i.e., homogeneous sets are uniformly thick with respect to
Lebesgue measure.
For instance, 
Cantor sets of positive length are homogeneous (see {\it Example} in
Sect. 1).  

The following space of character--automorphic forms (see e.g. [36])
will play an important role.

\proclaim{Definition} Let $E$ be a homogeneous set and let
$z:\Bbb D/\Gamma\equiv\Omega$ be a uniformization of the domain
$\Omega=\bar\Bbb C\setminus E$.
The space $A^{2}_1(\Gamma,\alpha)$ is
formed by functions $f$, which are analytic on $\Bbb D$ and satisfy the following
three conditions
$$
\align
1)& f \ \text{is of Smirnov class}\\
2)& \frac{f(\gamma(\zeta))}{\gamma_{21}\zeta+\gamma_{22}}=\alpha(\gamma)
f(\zeta)\quad
\forall
\gamma=\bmatrix \gamma_{11}&\gamma_{12}\\ \gamma_{21}&\gamma_{22}
\endbmatrix\in \Gamma\\ 3)& 
\int_{\Bbb E}\vert f\vert^{2}\,dm<\infty,
\endalign
$$
where
 $\Bbb E$ is a
fundamental set for the action of $\Gamma$ on $\Bbb T$.
\endproclaim

Let us mention that in our context the space $A_1^2(\Gamma,\alpha)$ arises
naturally in the following way. Since $\rho(t)$ from (0.3)
is supposed to satisfy Szeg\"o's
condition (0.7) it
 can be represented  in the form 
 $$
 \rho(t)=\vert D(t)\vert^2,
 $$
where $D(\zeta)$ is an outer function. Since the measure $\rho(t)\,dm(t)$
is invariant with respect to the substitution $t\to\gamma(t)$, we have
$$
\left\vert\frac{D(\gamma(t))}{\gamma_{21}t+\gamma_{22}}\right\vert^2=
\vert D(t)\vert^2,
$$
and hence the outer function $D(\zeta)$ itself satisfies an automorphic property
of the form
$$
\frac{D(\gamma(t))}{\gamma_{21}t+\gamma_{22}}=
\beta(\gamma)D(t),\quad \beta\in\Gamma^*.\tag 0.15
$$

$A^2_1(\Gamma,\alpha)$ is a closed subspace of $L^2_{dm\vert\Bbb E}$
with the reproducing kernel
$k^\alpha(\zeta,\zeta_0)$
(the point evaluation functional is bounded):
$$
\langle f(t), k^\alpha(t,\zeta_0)\rangle=f(\zeta_0), \quad
\zeta_0\in\Bbb D,\ f\in
A^2_1(\Gamma,\alpha).
$$
Put 
$$
k^\alpha(\zeta)= k^\alpha(\zeta,0)\quad\text{ and}\quad
K^\alpha(\zeta)=\frac{k^\alpha(\zeta)}{\sqrt{k^\alpha(0)}}.
$$

\proclaim{Theorem [29]} Let $E$ be a homogeneous set.
Let $z:\Bbb D/\Gamma\equiv\bar \Bbb C\setminus E$ with the normalization
$ (b z)(0)=1$.
Then the
systems of functions $\{b^n K^{\alpha\mu^{-n}}\}_{n\in\Bbb Z_+}$ 
and $\{b^n K^{\alpha\mu^{-n}}\}_{n\in\Bbb Z}$ 
form an
orthonormal basis in $A^2_1(\Gamma,\alpha)$ and
 in $L^2_{dm\vert\Bbb E} $, respectively,  for any $\alpha\in\Gamma^*$.
With respect to this basis, the operator multiplication by $z(t)$ is a
three--diagonal almost periodic Jacobi matrix, moreover
$$
z b^n K^{\alpha\mu^{-n}}=\Cal P(\alpha\mu^{-n}) b^{n-1} K^{\alpha\mu^{-n+1}}+
\Cal Q(\alpha\mu^{-n}) b^{n} K^{\alpha\mu^{-n}}+\Cal P(\alpha\mu^{-n-1})
b^{n+1} K^{\alpha\mu^{-n-1}},
$$
where
$$
\Cal P(\alpha)=\left(\frac{K^\alpha}{K^{\alpha\mu}}\right)(0),
\quad
b'(0)\Cal Q(\alpha)=(zb)'(0)+
\left(\log\frac{K^\alpha}{K^{\alpha\mu}}\right)'(0).
$$

Conversely, every almost periodic Jacobi matrix (0.14) such that
$\sigma(J)=\sigma_{a.c.}(J)=E$ can be represented in the form
$$
p_n=\Cal P(\alpha^{-1}\mu^{n+1}),\quad q_{n-1}=\Cal
Q(\alpha^{-1}\mu^{n+1}),\tag 0.16
$$
with some $\alpha\in\Gamma^*$.
\endproclaim

Let us point out, that the recurrence coefficients 
$\{p_n\}$ and $\{q_n\}$ from (0.16)
are uniformly almost periodic,
since $\Cal P(\alpha)$ and $\Cal Q(\alpha)$ are continuous functions on
the compact Abelian group $\Gamma^*$ (see (0.12)).

\medskip
In this paper for homogeneous sets $E$ questions a)--d) will be answered. One
of the main outputs of this paper is the following corollary of our Main
Theorem. For simplicity here in the introduction we restrict ourselves to the
case when no point--measures appear outside of $E$.

\proclaim{Corollary 0.1} Let $\sigma$ be a positive measure, whose support is
a homogeneous set $E$. Assume that $\log\sigma'_{a.c.}(z(t))\in
L^1$ and define an outer function $D(\zeta)$, $D(0)>0$, by the relation
$$
\vert D(t)\vert^2=2\pi\sigma'_{a.c.}(z(t))\vert z'(t)\vert,\quad t\in\Bbb T.
\tag 0.17
$$
Then the minimum deviation and the orthonormal polynomials
$P_n(z,\sigma)=\frac{z^n}{p_0\dots p_n}+\dots$ have the
following asymptotic behavior ($n\to\infty$)
$$
p_0\dots p_n \sim\frac{D(0)}{K^{\beta\mu^{n}}(0)},\tag 0.18
$$
$$
 P_n(z(\zeta))\sim\frac{ b^{-n}(\zeta)K^{\beta\mu^{n}}(\zeta)}
{D(\zeta)}
\tag 0.19
$$
uniformly on each compact subset of $\bar\Bbb C\setminus E$ and
$$
 D(t) P_n(z(t))-\left\{(b^{-n}K^{\beta\mu^{n}})(t)+\frac{D(t)}{D(\bar t)}
(b^{-n}K^{\beta\mu^{n}})(\bar t)\right\}\to 0\quad
\text{in}\  L^2_{dm\vert\Bbb E},\tag 0.20
$$
where $\beta$ and $\mu$ are given in (0.15) and (0.7), respectively.

\endproclaim

For a  comparison of this result with those one of
Szeg\"o and Widom (see (0.6), (0.10), (0.11)) let us note, that
for a homogeneous set $E$ we have the following analogue of (0.4):
$$
\text{clos}_{L^2_{dm\vert\Bbb E}}\{D H^\infty(\Gamma,\alpha)\}=
A_1^2(\Gamma,\beta\alpha).
$$
Thus the extremal function $f(\zeta,\rho,\alpha)$ in (0.9),
$\rho$ given by (0.8), is of the form
$$
f(\zeta,\rho,\alpha)=\frac{K^{\alpha\beta}(\zeta)}{K^{\alpha\beta}(0)}
\frac{D(0)}{D(\zeta)},
$$
and hence
$$
\nu(\rho,\alpha)=\left[\frac{D(0)}{K^{\alpha\beta}(0)}\right]^2.
$$
But let us point out that neither the ideas of proof nor the
methods of proof are related to  Widom's paper [34], only the
results can be considered as an extension of those one of Widom.

For the very special case $E=[-2,2]$ and thus $z(\zeta)=\zeta+1/\zeta$
the Fuchsian group $\Gamma$ is trivial, i.e., $\Gamma=\{\text{id}\}$,
$\Gamma^*=\{1\}$, $b(\zeta)=\zeta$, $K^\alpha(\zeta)\equiv 1$. Hence 
the basis $\{b^n K^{\alpha\mu^{-n}}\}$ becomes the standard Fourier basis
$\{t^n\}$ in $L^2$ and  Corollary 0.1 becomes Szeg\"o's result (in
particular (0.19) becomes (0.6)).

What concerns (0.20), it reminds strongly to the well known formulas from 
scattering theory [19, 21]. Roughly speaking, the
proof of our result
 is based on asymptotic orthogonality of "incoming" and "outgoing"
subspaces. More precisely we prove a character--automorphic analog of the
following well known proposition:
$$
P_+\{f(t) t^{-n}\}\to 0,\quad n\to\infty,
$$
where $f\in L^\infty$ and $P_+$ is the Riesz projection from $L^2$
onto $H^2$ (see Lemma 5.3).

Finally we would like to mention that even if $\sigma$ is absolutely
continuous  on
$E$ and
$\sigma'$ satisfies the Szeg\"o--type condition (0.3), then a
one--dimensional perturbation of the corresponding Jacobi matrix may lead to a
measure with a denumerable set of mass--points  outside $E$. Therefore it is
more natural to consider  measures which may have  a denumerable
set of mass--points  outside
$E$, where the mass--points are supposed to 
 accumulate on $E$ only. In fact, this point of view gives
us more freedom
and even helps us  to prove our Main Theorem. If such  mass--points appear
then  we assume that the set of poles (i.e. the set of mass--points) and zeros
of the associated Stieltjes function have to satisfy the
free--interpolation Carleson--type condition. Under these conditions 
asymptotics with respect to this wide class of measures are given.
With the help of the new asymptotic representation
the limit almost periodic behavior of the recurrence coefficients is
proved. 

\medskip

The paper is organized as follows: First the necessary ingredients from the
theory
of Hardy spaces of character--automorphic functions and forms are given. In the
second Section properties and a special representation of the Stieltjes
function are presented. Using this representation in Section 3 we introduce an
analogue of a scattering function and a transformation which plays a central
role in the proof. It maps polynomials in character--automorphic forms. 
Further, the main ideas of the proof of the Main Theorem are briefly
outlined. In Sections 4 and 5 we show that for a wide class of measures
$\sigma$ this map is a bounded map from $L^2_{d\sigma}$ to
$A_1^2(\Gamma,\alpha)$.  Approximating a
given measure by a sequence of such measures 
we are able to prove  the Main Theorem
and its Corollaries in the last Section.

\bigskip
\noindent{\bf Acknowledgment.} We would like to thank V.A. Marchenko for helpful
and stimulating discussions and M.Sodin
for explanations of some properties
of homogeneous sets.


\head 1. Preliminaries: the Hardy spaces on a Riemann surface of Widom
type
\endhead

For a compact set $E\subset\Bbb R$ let us consider
the open unit disk $\Bbb D$ as universal covering surface for the domain
$\Omega=\bar\Bbb C\setminus E$.
Thus there exists a meromorphic function $z(\zeta)$ mapping $\Bbb D$ onto
$\Omega$ and
a discrete subgroup $\Gamma$ of the group $SU(1,1)$ consisting of elements
of the form
$$
\gamma=\bmatrix \gamma_{11}&\gamma_{12}\\ \gamma_{21}&\gamma_{22}
\endbmatrix,\ 
\gamma_{11}=\overline{\gamma_{22}},\ \gamma_{12}=\overline{\gamma_{21}},
\ \det\gamma=1,
$$
such that 
 $z$ is automorphic with respect to $\Gamma$, i.e.,
$z(\gamma(\zeta))=z(\zeta),\ \forall \gamma\in\Gamma$,
and any two preimages of $z_0\in\Omega$ are $\Gamma$--equivalent, i.e.,
$$
z(\zeta_1)=z(\zeta_2)\ \Rightarrow\ \exists\gamma\in\Gamma:\ 
\zeta_1=\gamma(\zeta_2).
$$

As usual let us define 
$$
H^\infty(\Gamma)= \{f\in H^\infty:
\ f\circ\gamma=f,\ \forall \gamma\in\Gamma\}.
$$
Note, if the space $H^\infty(\Gamma)$ is not
trivial, i.e.,
$$
\exists f\in H^\infty(\Gamma):\ f(\zeta)\not\equiv f(\zeta_0),
$$
then the trajectory
$\{\gamma(\zeta_0)\}_{\gamma\in\Gamma}$ satisfies the Blaschke condition.
The Blaschke product
$$
 b(\zeta,\zeta_0)=b(\zeta,\zeta_0;\Gamma)=\prod_{\gamma\in\Gamma} 
 \frac{\gamma(\zeta_0)-\zeta}{1-\overline{\gamma(\zeta_0)}\zeta}
 \frac{\vert\gamma(\zeta_0)\vert}{\gamma(\zeta_0)}
$$
is called the {\it Green's function} of $\Gamma$ with respect to
$\zeta_0$.

If $G(z,z_0)$ denotes, as before, the Green's function of the domain
$\Omega$, then
$$
G(z(\zeta),z(\zeta_0))=-\log\vert b(\zeta,\zeta_0)\vert.
$$
The Green function is a character--automorphic function, that is there
exists a
$\mu_{\zeta_0}\in
\Gamma^*$ such that $ b(\gamma(\zeta),\zeta_0)=
\mu_{\zeta_0}(\gamma)b(\zeta,\zeta_0)$. To simplify the notation we put
$$
b(\zeta)=b(\zeta, 0)\quad\text{ and}\quad \mu=\mu_{0}.
$$

We will consider spaces of character--automorphic functions. For
$\alpha\in\Gamma^*$, define
$$
H^\infty(\Gamma,\alpha)= \{f\in H^\infty:
\ f\circ\gamma=\alpha(\gamma)f,\ \forall \gamma\in\Gamma\}.
$$
The domain $\bar\Bbb C\setminus E$ (respectively the group $\Gamma$) is said
to be of {\it Widom type} if for any
$\alpha\in\Gamma^*$ the space $H^\infty(\Gamma, \alpha)$ is not trivial, i.e.
$H^\infty(\Gamma, \alpha)\not=\{\text{const}\}$ [35, 28].
Also, in this case,  $\Gamma$ acts dissipative on $\Bbb T$
with respect to $dm$, that is there exists a measurable (fundamental)
set $\Bbb E$, which does not contain any two $\Gamma$--equivalent points,
and the union $\cup_{\gamma\in\Gamma}\gamma(\Bbb E)$ is a set of 
full measure [28].

Let $f$ be an analytic function in $\Bbb D$, $\gamma\in\Gamma$ and $k\in \Bbb N$.
Then we put
$$
f\vert[\gamma]_k=\frac{f(\gamma(\zeta))}{(\gamma_{21}\zeta+\gamma_{22})^k}
$$
It is easily verified that
$$
f\vert[\gamma_1\gamma_2]_k=(f\vert[\gamma_1]_k)\vert[\gamma_2]_k.
$$
Notice that $f\vert[\gamma]_2=f\ \forall\gamma\in \Gamma$, means that the form
$f(\zeta)d\zeta$ is invariant with respect to the substitutions
$\zeta\to\gamma(\zeta)$ ($f(\zeta)d\zeta$ is an Abelian
integral on $\Bbb D/\Gamma$). Analogically to $A_1^2(\Gamma,\alpha)$
(see Introduction) for a group of Widom type we define the space
$A_2^1(\Gamma,\alpha)$.

\proclaim{Definition} Let $\Gamma$ be a group of
Widom type and let $\Bbb E\subset \Bbb T$ be a fundamental set.
The space $A^{1}_2(\Gamma,\alpha)$
is formed by functions $f$, which are analytic on $\Bbb D$ and satisfy
the following three conditions
$$
\align
1)& f \ \text{is of Smirnov class}\\
2)& f\vert[\gamma]_2=\alpha(\gamma) f\quad \forall\gamma\in \Gamma\\
3)& 
\int_{\Bbb E}\vert f\vert\,dm<\infty.
\endalign
$$
\endproclaim

If $\Gamma$ is a group of Widom type then $A^2_1(\Gamma,\alpha)$ is a Hilbert
space with the reproducing kernel
$k^\alpha(\zeta,\zeta_0)$, moreover
$$
0<\inf_{\alpha\in\Gamma^*} k^\alpha(\zeta_0,\zeta_0)\le
\sup_{\alpha\in\Gamma^*} k^\alpha(\zeta_0,\zeta_0)<\infty.
$$
Let us also mention that the identity 
$$
\text{clos}_{L^2_{dm\vert\Bbb E}}\{D(t)H^\infty(\Gamma)\}=A_1^2(\Gamma,\alpha),
$$
where $D(t)$ is an outer function from $A_1^2(\Gamma,\alpha)$, generally speaking
does not hold for an arbitrary group of this type. It is valid for
groups of Widom type with Direct Cauchy Theorem.

\proclaim{Theorem [14]}
Let $\Gamma$ be a group of Widom type. The following statements are equivalent:

\noindent 1) The function $K^\alpha(0)$ is continuous on $\Gamma^*$.

\noindent 2) For $\alpha\in \Gamma^*$, let $\Delta^\alpha(\zeta)\in
H^\infty(\Gamma,\alpha)$, $\Vert\Delta^\alpha\Vert\le 1$,
be an extremal function of the problem
$$
\Delta^\alpha(0)=\sup\{\vert f(0)\vert:
f\in
H^\infty(\Gamma,\alpha),\ \Vert f\Vert\le 1\}.
$$
Then $\Delta^\alpha(0)\to 1$ ($\alpha\to 1_{\Gamma^*}$). 

\noindent 3)
The {\it Direct Cauchy Theorem} holds:
$$
\int_{\Bbb E}\frac{f}{b}(t)\,\frac{d t}{2\pi
i} =\frac{f}{b'}(0),
\quad\forall f\in A^1_2(\Gamma,\mu).\tag DCT
$$

\noindent 4) Let $\overline{t  A^2_1(\Gamma,\alpha^{-1})}=
 \{g=\overline{t f}:\ f\in A^2_1(\Gamma,\alpha^{-1})\}$. Then
$$
L^2_{dm\vert\Bbb E} =\overline{t  A^2_1(\Gamma,\alpha^{-1})}\oplus
 A^2_1(\Gamma,\alpha)\quad \forall\alpha\in\Gamma^*.
 $$
 
 \noindent 5) Every invariant subspace  $M\subset A^2_1(\Gamma,\alpha)$
(i.e. $f M\subset M\ \forall f\in H^\infty(\Gamma)$)
is of the form
$$ M=s A^2_1(\Gamma,\sigma^{-1}\alpha)$$ for some character--automorphic
inner function $s\in H^\infty(\sigma)$.

\endproclaim

\proclaim{Definition [6]}
A measurable set $E$ is homogeneous if there is an $\eta>0$ such that
$$
\vert(x-\rho,x+\rho)\cap E\vert\ge \eta\rho\quad \text{for all}\
0<\rho<\text{\rm diam}\, E \quad\text {and all}\ x\in E.\tag 1.1
$$
\endproclaim

\demo{Example} Let us demonstrate that Cantor sets of positive length
are homogeneous. First recall the
construction of such a set [20]. We start with an
interval of the length
$l_0$ and take a sequence of numbers
$$
\{\varkappa_j\}_{j\ge 1}:\quad 0<\varkappa_j<1,\ 
\sum_{j\ge1}\varkappa_j<\infty.
$$
At the first step we remove an open segment,
whose length is $\varkappa_1$ part of the common length of the initial
interval, so that on either side there remains a closed line segment
of the length $l_1=\frac 1 2{(1-\varkappa_1)l_0}$. Then we make the same
procedure with each of the remaining intervals, taking out $\varkappa_2$
part of each of them. Continuing in this way, we get a Cantor set
$E=E(l_0;\varkappa_1,\varkappa_2,\dots)$, with the common length
$$
\vert E(l_0;\varkappa_1,\varkappa_2,\dots)\vert=
\prod_{j\ge 1}(1-\varkappa_j) l_0>0.
$$
Let us mention that this set consists of two Cantor sets of the
form $E(l_1;\varkappa_2\dots)$, with $l_1=\frac 1 2{(1-\varkappa_1)l_0}$,
or of four sets of the form
$E(l_2;\varkappa_3,\dots)$, $l_2=\frac 1 2{(1-\varkappa_2)l_1}$,
and so on, ..., of  $2^n$ sets of the form
$E(l_n;\varkappa_{n+1},\dots)$, $l_n=\frac 1 2{(1-\varkappa_n)l_{n-1}}$.

Now we check that
$$
\vert(x-\rho,x+\rho)\cap E\vert\ge\frac 
{\prod_{j\ge 1}(1-\varkappa_j)} 2\rho,\quad \forall x\in E,\ 
\forall\rho\le
\text{diam}\, E.
$$
Let $\frac 1 2{\text{diam}\, E} \le\rho\le\text{diam}\, E$, and
$x\in E$. Then the interval $(x-\rho,x+\rho)$ contains at least
one of the set $E(l_1;\varkappa_2,\dots)$. So,
$$
\align
\vert(x-\rho,x+\rho)\cap E\vert\ge&\vert E(l_1;\varkappa_2,\dots)\vert=
\prod_{j\ge 2}(1-\varkappa_j)\cdot(1-\varkappa_1)\frac{l_0}2\\
=&\frac{\prod_{j\ge 1}(1-\varkappa_j)}2 l_0\ge
\frac{\prod_{j\ge 1}(1-\varkappa_j)}2 \rho.
\endalign
$$
If $\frac 1 4{\text{diam}\, E} \le\rho\le\frac 1 2\text{diam}\, E$
and $x\in E$, then the interval $(x-\rho,x+\rho)$ contains at least
one of the set $E(l_2;\varkappa_3,\dots)$. Finally,
if $\frac 1 {2^n}{\text{diam}\, E} \le\rho\le\frac 1
{2^{n-1}}\text{diam}\, E$ and $x\in E$, then the interval
$(x-\rho,x+\rho)$ contains at least one of the set
$E(l_n;\varkappa_{n+1},\dots)$. Hence,
$$
\vert(x-\rho,x+\rho)\cap E\vert\ge\vert E(l_n;\varkappa_{n+1},\dots)\vert
=\frac{\prod_{j\ge 1}(1-\varkappa_j)}{2^n} l_0\ge
\frac{\prod_{j\ge 1}(1-\varkappa_j)}2 \rho.
$$
\qed
\enddemo

The following sufficient condition of homogeneity was proposed and proved
by M. Sodin based on [8].

\proclaim{Proposition (M. Sodin)} Let $E=[b_0,a_0]\setminus\cup_{j\ge 1}
(a_j,b_j)$. Let $l_j=b_j-a_j$, $j\ge 1$ and put, formally,
$l_0=1$. Let $\rho_{j,k}$ be the distance between the two gaps $(a_j,b_j)$
and $(a_k,b_k)$, $j\not=k,\ j,k\ge 1$ and let $\rho_{j,0}$ denote the distance
from the gap $(a_j,b_j)$ to the boundary of the interval
$[b_0,a_0]$. If
$$
\sup_j\sum_{k\not=j}\frac{l_j^{1/2}l_k^{1/2}}{\rho_{j,k}}<\infty,
$$
then $E$ is a homogeneous set.
\endproclaim

Let us mention that the case of Julia sets of polynomials
which are
real and thus (up to the classical Chebyshev polynomials and its
conjugates) of Cantor type are not homogeneous because the length of the
Julia set is zero.

\proclaim{Theorem} Let $E$ be a homogeneous set, 
then $\bar\Bbb C\setminus E$ is
of Widom type and  the Direct Cauchy Theorem holds.
\endproclaim

The proof of this theorem is based mainly on the following
lemma of Jones and Marshall [16], who following Carleson,
considered the Corona problem
for the surface $\bar\Bbb C\setminus E$,
where $E$ is a homogeneous set.

\proclaim{Lemma} Let $E=[b_0,a_0]\setminus\cup_{j\ge 1}(a_j,b_j)$
be a homogeneous set. From each interval $(a_j,b_j)$ let us pick arbitrarily 
exactly one point
$x_j$. Then there is a constant $N$
which depends only on the value $\eta$ in (1.1), such
that
$$
\sup_{j\ge 1}\sum_{i\not= j}G(x_i,x_j)\le N<\infty.
$$
\endproclaim

\bigskip

To finish this section,
we would like to extend the DCT to  functions  with infinitely
many  poles.

\proclaim{Lemma 1.1} Let $B\in H^\infty(\Gamma,\alpha)$ be a Blaschke
product,
$B=\prod_l b(\zeta,\zeta_l)$,  and
$f\in A_2^1(\Gamma,\alpha)$. If
$$
\sum_l\left\vert\frac{f(\zeta_l)}{B'(\zeta_l)}\right\vert<\infty,
$$
then
$$
\int_{\Bbb E}\frac{f}{B}\frac{dt}{2\pi i}=
\sum_l\frac{f(\zeta_l)}{B'(\zeta_l)}.
$$
\endproclaim

\demo{Proof} In the proof we use a Poisson--like kernel. Let
$K^\beta(\zeta,\zeta_l)=\frac{k^\beta(\zeta,\zeta_l)}
{\sqrt{k^\beta(\zeta_l,\zeta_l)}}$. Then
$$
1=\int_{\Bbb E}\vert K^\beta(t,\zeta_l)\vert^2\,dm=\int_{\Bbb
E} K^\beta(t,\zeta_l) (\overline{t K^\beta(t,\zeta_l)})\frac{dt}{2\pi i}.
$$
Let us show that 
$$
\overline{t
K^\beta(t,\zeta_l)}=
\frac{K^{\beta^{-1}\mu_{\zeta_l}}(t,\zeta_l)}{b(t,\zeta_l)}
\frac{b'(\zeta_l,\zeta_l)}{\vert b'(\zeta_l,\zeta_l)\vert}.
\tag 1.2
$$
The function 
$ \frac{K^{\beta^{-1}\mu_{\zeta_l}}(\zeta,\zeta_l)}{b(\zeta,\zeta_l)}$ is
orthogonal to $A_1^2(\Gamma,\beta^{-1})$, and, hence, it is of the form
$ \frac{K^{\beta^{-1}\mu_{\zeta_l}}(t,\zeta_l)}{b(t,\zeta_l)}=
\overline{tg_0}$, $g_0\in A_1^2(\Gamma,\beta)$. But, due to DCT,
$$\align
\langle g,g_0\rangle=&\langle\overline{t g_0},\overline{t g}\rangle=
\left\langle\frac{K^{\beta^{-1}\mu_{\zeta_l}}(t,\zeta_l)}{b(t,\zeta_l)},
\overline{t g}\right\rangle\\=&
\int_{\Bbb
E} \frac{K^{\beta^{-1}\mu_{\zeta_l}}(t,\zeta_l)}{b(t,\zeta_l)}
g(t)\frac{dt}{2\pi i}=
\frac{K^{\beta^{-1}\mu_{\zeta_l}}(\zeta_l,\zeta_l)}{b'(\zeta_l,\zeta_l)}
g(\zeta_l).
\endalign
$$
Therefore,
$$
g_0=k^\beta(t,\zeta_l)\overline{\left(
\frac{K^{\beta^{-1}\mu_{\zeta_l}}(\zeta_l,\zeta_l)}{b'(\zeta_l,\zeta_l)}
\right)}
=K^\beta(t,\zeta_l)
\frac{K^\beta(\zeta_l,\zeta_l)K^{\beta^{-1}\mu_{\zeta_l}}(\zeta_l,\zeta_l)}
{\vert b'(\zeta_l,\zeta_l)\vert}
\frac{b'(\zeta_l,\zeta_l)}{\vert b'(\zeta_l,\zeta_l)\vert}.
$$
And, since $\Vert g_0\Vert=1$, we get (1.2), and the identity
$$
\frac{K^\beta(\zeta_l,\zeta_l)K^{\beta^{-1}\mu_{\zeta_l}}(\zeta_l,\zeta_l)}
{\vert b'(\zeta_l,\zeta_l)\vert}=1.
$$
Put
$$
P^\beta(t,\zeta_l)=\frac{K^\beta(t,\zeta_l)
K^{\beta^{-1}\mu_{\zeta_l}}(t,\zeta_l)}{b(t,\zeta_l)}
\frac{b'(\zeta_l,\zeta_l)}{\vert b'(\zeta_l,\zeta_l)\vert}.
$$
As it was shown $P^\beta(t,\zeta_l)\frac{dt}{2\pi i}\ge 0$, and
$\int_{\Bbb E}P^\beta(t,\zeta_l)\frac{dt}{2\pi i}= 1$.

\medskip

For the given functions $B$ and $f$ consider the series
$$
\tilde f=\sum_l\frac{f(\zeta_l)}{B'(\zeta_l)}P^\beta(t,\zeta_l)B,
$$
with some $\beta\in\Gamma^*$. This series converges absolutely
in $A^1_2(\Gamma,\alpha)$ and interpolates the function $f$ since
$$
\tilde f(\zeta_l)=\frac{f(\zeta_l)}{B'(\zeta_l)}
\frac{K^\beta(\zeta_l,\zeta_l)K^{\beta^{-1}\mu_{\zeta_l}}(\zeta_l,\zeta_l)}
{ b'(\zeta_l,\zeta_l)}
\frac{b'(\zeta_l,\zeta_l)}{\vert b'(\zeta_l,\zeta_l)\vert}
B'(\zeta_l)=f(\zeta_l).
$$
So, the function $f$ can be represented in the  form:
$f=\tilde f+B g$, with $g\in A_2^1(\Gamma)$. Due to DCT
$$
\int_{\Bbb E}\frac{f}{B}\frac{dt}{2\pi i}-\int_{\Bbb E}\frac{\tilde
f}{B}\frac{dt}{2\pi i}=
\int_{\Bbb E}\frac{f-\tilde f}{B}\frac{dt}{2\pi i}=0,
$$
and since the series converges absolutely we can integrate it term by
term, which proves the lemma.\qed
\enddemo


\head{2. A special representation of the Stieltjes function}
\endhead

For the following we need some additional notation.
Let $E$ be a homogeneous set and $X\subset \Bbb R\setminus E$ be
a set of points which can accumulate only to the set $E$. Let 
$\sigma$ be a positive measure with  support $E\cup X$. To this measure we
associate the so--called Stieltjes function with a special normalization
$$
r(z)=1+\int\frac{d\sigma(x)}{x-z}
=1+\sum_{x_l\in X}\frac{\sigma_l}{x_l-z}+
\int_E\frac{d\sigma(x)}{x-z}.\tag 2.1
$$
It is a function meromorphic on $\Omega=\bar\Bbb C\setminus E$ and such that
$$
\frac{r(z)-\overline{r(z)}}{z-\bar z}\ge 0.
$$
The function $-1/r(z)$ possesses the same properties, and therefore it has the
representation
$$
-1/r(z)=-1+\int\frac{d\sigma^{(\tau)}(x)}{x-z}
=-1+\sum_{x_l^{(\tau)}\in X^{(\tau)}}\frac{\sigma_l^{(\tau)}}{x_l^{(\tau)}-z}+
\int_E\frac{d\sigma^{(\tau)}(x)}{x-z}
.\tag 2.2
$$
The support of the measure $d\sigma^{(\tau)}(x)$ is the set
$E\cup X^{(\tau)}$, where $X^{(\tau)}$ is the set of zeros of the function
$r(z)$.

The following property of the measures $\sigma$ and $\sigma^{(\tau)}$
plays an essential role in what follows.

\proclaim{Lemma 2.1} Suppose that the measures $\sigma$ and $\sigma^{(\tau)}$ are
related by (2.1) and (2.2). Then the polynomial map
$$
P^{(\tau)}(z)=P(z)-\int\frac{P(x)-P(z)}{x-z}\,d\sigma(x)\tag 2.3
$$
induces a unitary map from
$L^2_{d\sigma}$ to $L^2_{d\sigma^{(\tau)}}$. Furthermore, the inverse operator
has the form
$$
P(z)=P^{(\tau)}(z)+\int\frac{P^{(\tau)}(x)-P^{(\tau)}(z)}{x-z}\,d\sigma^{(\tau)}(x).\tag
2.4
$$
\endproclaim

\demo{Proof} Let $J$ be the Jacobi matrix associated with the measure $\sigma$,
$$
\int\frac{d\sigma(x)}{x-z}=\langle (J-z)^{-1} e,e\rangle,
$$
where
$$
J=\bmatrix q_0&p_1&0&0\\
p_1&q_1&p_2&0\\
0&p_2&q_3&\ddots\\
0&0&\ddots&\ddots
\endbmatrix
\quad\text{and}\quad
e=\bmatrix p_0\\ 0\\0\\ \vdots\endbmatrix.
$$
Direct calculation shows that
$$
\int\frac{d\sigma^{(\tau)}(x)}{x-z}=\langle (J^{(\tau)}-z)^{-1} e,e\rangle,
$$
where
$$
J^{(\tau)}= J+e\langle\ ,e\rangle.
$$
Therefore, orthonormal polynomials with respect to the measure $\sigma^{(\tau)}$
satisfy the same recurrence relation (0.2) but with
initial data
$$
p_0 P^{(\tau)}_0(z)=1,\quad zP_0^{(\tau)}=(q_0+p_0^2)P_0^{(\tau)}+P^{(\tau)}_1.
$$
For this reason (2.3) transforms the system of orthonormal
polynomials corresponding to the measure $\sigma$
into the system
of orthonormal polynomials with respect to the measure $\sigma^{(\tau)}$.

The same arguments show that (2.4) is the inverse map.\qed

\enddemo

The following theorem describes quite general properties
of a function $r(z(\zeta))$ as a function on the universal covering.

\proclaim{Theorem [29]} Assume that the set $X$ of poles of $r(z)$ 
satisfies the condition
$$
\sum_{x_l\in X} G(x_l)<\infty.\tag 2.5
$$
then $r(z(\zeta))$ is a function of
bounded characteristic without a singular component.
\endproclaim

We will use essentially a special representation of $r(z(\zeta))$,
which follows from this Theorem.

\bigskip

\proclaim{Lemma 2.2} Let $r(z)$ be a function  of the form
(2.1).
Suppose that $X$ satisfies
 condition (2.5) and assume that $\log\sigma'_{a.c.}(z(t))\in L^1$.

Then $r(z)$ has a  representation of the form
$$
r(z(\zeta))=\frac{\psi(\zeta)}{\phi(\zeta)},\tag 2.6
$$
where $\psi$ and $\phi$ are of Smirnov class;
$\psi\vert[\gamma]=\alpha(\gamma)\psi$, $\phi\vert[\gamma]=\alpha(\gamma)\phi$
with some $\alpha\in\Gamma^*$ and
$$
\psi(t)\overline{\phi(t)}-\phi(t)\overline{\psi(t)}=-tz'(t),\quad
t\in\Bbb T.\tag 2.7
$$
\endproclaim

\demo{Proof} 
We would like to define an outer part of $\phi$ by the relation
$$
\vert\phi\vert^2\text{Im}\ r\circ z=-\frac{t z'}{2i}.\tag 2.8
$$
It is known that $z'$ is of bounded characteristic, moreover
$b^2 z'$ is the outer function.
Since $\frac 1 \pi \text{Im}\ r=\sigma'_{a.c.}$ and $\log\sigma'_{a.c.}(z(t))\in
L^1$, an outer function $\phi^{\text{out}}$
is well defined by (2.8).

Since $\text{Im}\ r\circ z$ is an automorphic function, and since
$$
\left\vert \frac{z'\circ\gamma}{(\gamma_{21}t+\gamma_{22})^2}\right\vert=
\vert z'\circ\gamma\vert\vert\gamma'\vert=\vert z'\vert,
$$
we have
$$
\left\vert
\frac{\phi^{\text{out}}\circ\gamma}{\gamma_{21}t+\gamma_{22}}\right\vert^2
=\vert\phi^{\text{out}}\vert^2.
$$
Due to the uniqueness property of an outer function
$$
\frac{\phi^{\text{out}}\circ\gamma}{\gamma_{21}t+\gamma_{22}}
=\alpha_{\text{out}}(\gamma)\phi^{\text{out}},
$$
with some $\alpha_{\text{out}}\in\Gamma^*$.
Let $B_{X}(\zeta)$ denote the Blaschke product
 with zeros at $\{\zeta_l:\ z(\zeta_l)\in X\}$. Due to condition
(2.5) it is well defined. Now we put
$\phi(\zeta):=B_{X}(\zeta)\phi^{\text{out}}(\zeta)$. This is a function of
bounded characteristic without a singular component, possessing the automorphic
property we need. Next, let us put $\psi=(r\circ z)\phi$. This function has no
poles, and, due to the cited Theorem, it also does not have a singular
component. Thus the lemma is proved.\qed
\enddemo
\proclaim{Remark} We would like to point out that the function $\phi$ depends
only on the absolutely continuous part of the measure $\sigma_{a.c.}$ on $E$ and
on the support of the measure $\sigma$ outside $E$, i.e. on the set
$X$. We will see that all asymptotics are given only in terms of these
functions.

By the way, the outer function $D$
in (0.17) and the function $\phi$ are related in the following way
$\vert D\vert=\vert z'/\phi\vert$. So, under the normalization
$D(0)>0$ and $\phi^{\text{out}}(0)>0$, we have
$$
\frac 1 {D(\zeta)}=\left(\frac{\phi^{\text{out}}}{(-z')b^2}\right)(\zeta).
\tag 2.9
$$
\endproclaim
\bigskip

Besides the additive representation (2.1) the function $r(z)$ possesses 
the following
exponential (or multiplicative) representation.

\proclaim{Lemma 2.3} Let $r(z)$ be a function of the form (2.1).
Denote by $x_k^{(\tau)}\in X^{(\tau)}$ the nearest right-hand side 
zero to the pole $x_k\in X$.
Then
$$
r(z)=\exp\left\{\frac 1\pi\int\frac{f(x)\,dx}{x-z}\right\},\tag 2.10
$$
where
$$
f(x)=\left\{\aligned \pi\ \ ,&\ x_k<x<x_k^{(\tau)}\\
                     u(x),&\  x\in E\\
                     0\ \ ,&\ \text{ otherwise}\endaligned
\right. $$
and $u(x)=\arg r(x)$, $x\in E$.

Moreover, if (2.5) holds, then $\log\sigma'_{a.c.}(z(t))\in L^1$ if and only if
$\log\sin u (z(t))\in L^1$.
\endproclaim
\demo{Proof} The representation (2.10) follows immediately from
the general exponential representation of  functions
with positive imaginary part in the upper half-plane:
$$
r(z)=\text{Const}\cdot\exp\left\{\frac 1\pi\int_{\Bbb
R}\arg r(x+i0)\left(\frac{1}{x-z}-
\frac{x}{1+x^2}\right)\,dx\right\}.
$$
We only have to mention, that outside $E$ the function $r(z)$
is real, so the argument of $r(z)$ here is equal to $0$
or $\pi$ and note that $r(\infty)=1$.

 Splitting the integral
into two parts we have
$$\align
r(z)=&\exp\left\{\sum\int_{x_k}^{x^{(\tau)}_k}\frac{dx}{x-z}
\right\}\cdot\exp\left\{\frac 1\pi\int_{E}u(x)\frac{dx}{x-z}\right\}\\
=& r_1(z)\cdot r_2(z).
\endalign
$$
The first factor has real boundary values on $E$, moreover they are
positive, therefore $\text{Im}\ r(x+i0)=r_1(x+i0)\cdot\text{Im}\ r_2(x+i0)$.
Due to the Theorem mentioned before  $r_1\circ z$ is of bounded
characteristic, so $\log \vert r_1\circ z\vert=\log r_1\circ z\in L^1$.

Next let us consider the second factor. The function $\log r_2(z)$
is holomorphic in $\bar\Bbb C\setminus E$, its imaginary part
$$
U(z)=\frac 1\pi\int_{E}\frac{\text{Im}\ z}{\vert x-z\vert^2}u(x)\,dx
$$
is a bounded harmonic function. Consequently, $U\circ z$ is a bounded harmonic
function in $\Bbb D$. Therefore the boundary values of the conjugated
function $\log\vert r_2\circ z\vert$, for sure, belong to
$L^1$. Now, $\text{Im}\ r_2= \exp\{\log\vert r_2\vert\}
\sin\{\arg r_2\}= \exp\{\log\vert r_2\vert\}
\sin\{\arg r\}$. 
So $\log\sigma'_{a.c}\circ z=\log\{\frac 1\pi\text{Im}\ r\circ z\}$ is in
$L^1$ if and only if
$\log\sin\{u\circ z\}=\log\sin\{\arg r\circ z\}$ is in $L^1$.
\qed
\enddemo

The last lemma of this section gives us a sufficient condition such that the
measure $\sigma$ has no singular component on $E$.
\proclaim{Lemma 2.4} Assume that (2.5) holds and that
$$
\int_E\vert r(x)\vert\,dx<\infty.\tag 2.11
$$
Then $\sigma_{s.}\vert E=0$.
\endproclaim
\demo{Proof} Series expansion of $r(z)$ at infinity gives by (2.1)
$$
r(z)=1-\frac{p_0^2} z+\dots,\tag 2.12
$$
where
$$
p_0^2=\sum_l \sigma_l+\int_E\,d\sigma(x).
$$
It suffices to show that
$$
p_0^2=\sum_l \sigma_l+\frac 1\pi\int_E\text{Im}\ r(x+i0)\,dx.\tag 2.13
$$

Since $B_X r\circ z \cdot b^2 z'$ is a function
of Smirnov class, and due to (2.11) it is integrable we 
conclude that $B_X r\circ z \cdot b^2 z'\in A_2^1(\Gamma,\beta)$ with
some $\beta\in\Gamma^*$. The differential $r(z)\,dz$ has poles at the points
$x_l$ and at infinity. Taking in mind (2.1) and (2.12) we get, that the
sum of residues of the given differential is equal to
$$
\sum\text{Res}\ r(z)\,dz=-\sum_l\sigma_l+p_0^2.
$$
Since this series converges absolutely, using Lemma 1.1, we obtain
$$
\frac 1{2\pi i}\int_{\partial\Omega}r(z)\,dz=-\sum_l\sigma_l+p_0^2.
$$
Due to the symmetry property $\overline{r(\bar z)}=r(z)$, we have
$\text{Re}\ r(x+i0)=\text{Re}\ r(x-i0)$
and $\text{Im}\ r(x+i0)=-\text{Im}\ r(x-i0)$, hence
$$
\frac 1{2\pi i}\int_{\partial\Omega}r(z)\,dz=
\frac 1{\pi }\int_{E}\text{Im}\ r(x+i0)\,dx.
$$
Thus, (2.13) is proved.\qed
\enddemo

\head{3. From spectral data to scattering data}
\endhead

In this Section we will first derive some properties of the following map
$$
h(\zeta)=h(\zeta, P;\sigma)=\phi(\zeta)\int\frac{P(x)}{z(\zeta)-x}\,d\sigma(x),
\tag 3.1
$$
which will play an important role in what follows. Then we will briefly outline
the main ideas of the proof of the Main Theorem given in Section 6.

The map $h$ from (3.1) can be represented  in the form
$$\align
h(\zeta,P;\sigma)=&\phi(\zeta)
\left\{\int\frac{P(x)-P(z(\zeta))}{z(\zeta)-x}\,d\sigma(x)
+P(z(\zeta))\int\frac{d\sigma(x)}{z(\zeta)-x}\right\}\\
=&\phi(\zeta)
\left\{P(z(\zeta))-\int
\frac{P(x)-P(z(\zeta))}{x-z(\zeta)}\,d\sigma(x)-
P(z(\zeta))r(z(\zeta))\right\}\\
=&
\phi(\zeta) P^{(\tau)}(z(\zeta))-\psi(\zeta) P(z(\zeta)),\tag 3.2\endalign
$$
as well as by Lemma 2.1
$$
h(\zeta)=h(\zeta,
P;\sigma)=\psi(\zeta)\int\frac{P^{(\tau)}(x)}{z(\zeta)-x}\,d\sigma^{(\tau)}(x).
$$

With the help of these representations we obtain the following properties of
the function $h$:
$h\vert[\gamma]=\alpha(\gamma)h$, since the functions $\phi$ and $\psi$ have
this property; h is of Smirnov class. Indeed, if $n=\deg P$, then $b^n P\in
H^\infty$. So $b^n h(\zeta,P;\sigma)$ is of Smirnov class. It remains to be
shown, that
$h$ has no poles at $\{\gamma(0)\}$. In fact, it has a zero at the origin.
This
follows immediately from (3.1). Since we can not guarantee, that $h$ is
square-integrable on $\Bbb E$, in general $h $ does not belong to
$A_1^2(\Gamma,\alpha)$.

Using transformation (3.1) we will pass from the standard spectral parameter,
i.e., spectral measure $\sigma$, to a system of objects which remind to
scattering data, i.e., to a unimodular function $s$ and a discrete measure $\nu$
supported outside of the essential spectrum
(for the connection between scattering theory and orthogonal polynomials see
[21]). What we are going to do in fact, it's only to rewrite a standard norm
$\int\vert P\vert^2\,d\sigma$ in terms of the function $h$. Doing this at the
last step we assume that
$\sigma$ has no singular spectrum on $E$.

First let us express $P$ and $P^{(\tau)}$  in
terms of the function
$h$. Since
$\phi(\bar t)=\overline{\phi(t)}$ and
$\psi(\bar t)=\overline{\psi(t)}$, we get from  (3.2):
$$
\bmatrix h(t)\\h(\bar t)\endbmatrix=
\bmatrix-\psi&\phi\\ -\bar\psi&\bar\phi\endbmatrix
\bmatrix P\\ P^{(\tau)}\endbmatrix.
$$
Therefore,  using (2.6) we get
$$
tz'\bmatrix P\\ P^{(\tau)}\endbmatrix=
\bmatrix\bar\phi&-\phi\\ \bar\psi&-\psi\endbmatrix
\bmatrix h(t)\\h(\bar t)\endbmatrix,
$$
which gives the desired representation
$$
\bmatrix\overline{\left(\frac{t z'}{\phi}\right)} P\\ 
\overline{\left(\frac{t z'}{\psi}\right)}P^{(\tau)}\endbmatrix=-
\bmatrix 1& -s\\ 1& -s^{(\tau)}\endbmatrix
\bmatrix h(t)\\\bar t h(\bar t)\endbmatrix,\tag 3.3
$$
where $s=\phi/(\overline{t\phi})$,
$s^{(\tau)}=\psi/(\overline{t\psi})$.

 Since
$\sigma'_{a.c.}(x)=\frac 1\pi\text{Im}\ r(x+i0)$ it follows by (2.7),
$$
\aligned
\int_E \vert P(x)\vert^2\,d\sigma_{a.c.}(x)=
&\frac 1 2 \int_{\Bbb
E}\vert P(z(t))\vert^2\frac 1
\pi
\left(\frac{-tz'}{2 i\vert\phi\vert^2}\right)(z'\, dt)\\
=&\frac 1 2 \int_{\Bbb E}\vert P(z(t))\vert^2
\frac{(-tz')(t z')}{\vert\phi\vert^2}\, dm,
\endaligned\tag 3.4
$$
where $P$ is a polynomial.
Substituting (3.3) in (3.4) we get for the
absolutely continuous part of the measure $\sigma(x)$:
$$
\int_E \vert P(x)\vert^2\,d\sigma_{a.c.}(x)=\int_{\Bbb E}\left\vert
\frac{h(t)-s(t)\bar t h(\bar t)}{\sqrt 2}\right\vert^2\,dm(t).\tag 3.5
$$
Analogically, for the absolutely continuous part of the measure
$d\sigma^{(\tau)}$, using (2.6), we get
$$
\int_E \vert P^{(\tau)}(x)\vert^2\,d\sigma^{(\tau)}_{a.c.}(x)=
\int_{\Bbb E}\left\vert
\frac{h(t)-s^{(\tau)}(t)\bar t h(\bar t)}{\sqrt 2}\right\vert^2\,dm(t).\tag 3.6
$$

For the pure point spectrum of $\sigma$ and $\sigma^{(\tau)}$ we also can pass
from a polynomial
$P(z)$ to the function
$h(\zeta)=h(\zeta,P;\sigma)$.
Since $r=\psi/\phi$, we have $\sigma_k=-\left(\frac\psi{\phi'}z'\right)(x_k)$
and
$\sigma^{(\tau)}_k=\left(\frac\phi{\psi'}z'\right)(x^{(\tau)}_k)$.
In view of  (3.2),
$$
P(x_k)=-\frac h \psi(x_k)\quad\text{and}\quad
P^{(\tau)}(x^{(\tau)}_k)=\frac h \phi(x_k^{(\tau)}),
$$
so
$$
\sum_{x_k\in X}\vert P(x_k)\vert^2{\sigma_k}=\sum_{x_k\in X}\left\vert
\frac h \psi\right\vert^2(x_k){\sigma_k},
$$
and
$$
\sum_{x_k^{(\tau)}\in
X^{(\tau)}}\vert P^{(\tau)}(x_k^{(\tau)})\vert^2{\sigma_k^{(\tau)}}=
\sum_{x_k^{(\tau)}\in
X^{(\tau)}}\left\vert
\frac h \phi\right\vert^2(x_k^{(\tau)}){\sigma_k^{(\tau)}}.
$$
Thus we can define  discrete measures $\nu$ and $\nu^{(\tau)}$ on $\Bbb D$ 
$$
\nu(\zeta_l)=\frac{\sigma_l}{\vert\psi(\zeta_l)\vert^2},\quad
z(\zeta_l)=x_l,\tag 3.7
$$
and
$$
\nu^{(\tau)}(\zeta_l^{(\tau)})=\frac{\sigma_l^{(\tau)}}
{\vert\phi(\zeta_l^{(\tau)})\vert^2},\quad
z(\zeta_l^{(\tau)})=x_l^{(\tau)}.\tag 3.8
$$
These measures possess the
automorphic property 
$${\vert\gamma_{21}\zeta_l+\gamma_{22}\vert^2}
\nu(\gamma(\zeta_l))={\nu(\zeta_l)},
\quad
{\vert\gamma_{21}\zeta_l^{(\tau)}+\gamma_{22}\vert^2}
\nu^{(\tau)}(\gamma(\zeta_l^{(\tau)}))={\nu(\zeta_l^{(\tau)})}.\tag 3.9
$$
Let $Z$ be a fundamental set in the support of the measure $\nu$.
Then the sum
$$
\sum_{\zeta_l\in Z}\vert g(\zeta_l)\vert^2\nu(\zeta_l),\quad
g\in A_1^2(\Gamma,\beta),
$$
in fact, does not depend on the choice of the fundamental set and
$$
\sum_{x_k\in X}\vert P(x_k)\vert^2{\sigma_k}=\sum_{\zeta_l\in Z}
\vert h(\zeta_l)\vert^2\nu(\zeta_l).
$$
Analogically,
$$
\sum_{x_k^{(\tau)}\in X^{(\tau)}}\vert
P^{(\tau)}(x_k^{(\tau)})\vert^2{\sigma_k^{(\tau)}}=\sum_{\zeta_l^{(\tau)}
\in Z^{(\tau)}}
\vert h(\zeta_l^{(\tau)})\vert^2\nu^{(\tau)}(\zeta_l^{(\tau)}).
$$

Hence, under the assumption $\sigma_{s.}\vert E=0$, we get
$$
\int \vert P\vert^2\,d\sigma
=\int_{\Bbb E}\left\vert
\frac{h(t)-s(t)\bar t h(\bar t)}{\sqrt 2}\right\vert^2\,dm(t)
+\sum_{\zeta_l\in Z}
\vert h(\zeta_l)\vert^2\nu(\zeta_l)=:\Vert h\Vert_{s,\nu}^2.\tag 3.10
$$

\bigskip

Now let  us briefly and roughly outline
the proof of the Main Theorem. Assume that
$P(z)\to h(\zeta,P;\sigma)$ is a bounded map from $L^2_{d\sigma}$ to
$A_1^2(\Gamma,\alpha)$, i.e.,
$$
\int_{\Bbb E}\vert h\vert^2\,dm\le C_1\int\vert P\vert^2\,d\sigma\tag 3.11
$$
and that the measure $\nu$ possesses the property
$$
\sum_{\zeta_l\in Z}\vert h(\zeta_l)\vert^2\nu(\zeta_l)\le
C_2\int_{\Bbb E}\vert h\vert^2\,dm,\tag 3.12
$$
then 
$\Vert \cdot\Vert^2_{s,\nu}$ in (3.10)
gives us a norm in $A_1^2(\Gamma,\alpha)$ which is equivalent
to the original one,
$$
\frac 1{C_1}\Vert h\Vert^2\le\Vert h\Vert^2_{s,\nu}\le(2+C_2)\Vert h\Vert^2.
$$

\bigskip

Next let us consider the special
system of functions $\{h_n(\zeta)\}$,
$$
h_n(\zeta)=h(\zeta,P_n;\sigma),
$$
which is important in proving our results, as we shall see in a moment.
Note, that $h_n/\phi$ is the so called $n$--th function of second kind.
The
functions satisfy a three-term recurrence relation of the form (0.2),
$$
z(\zeta)h_n(\zeta)=p_n h_{n-1}(\zeta)+q_n h_{n}(\zeta)+p_{n+1}
h_{n+1}(\zeta),\tag 3.13
$$
and by orthogonality of $P_n$ (0.1) and
definition (3.1) it follows that $h_n(\zeta)$ has a zero of  order $n+1$ at
the origin. From this remark, and (3.13) we derive that 
$$
(z^{n+1}h_n)(0)=p_0\dots p_n\phi(0).\tag 3.14
$$

Thus, taking into consideration (3.10), we obtain
 that $\{h_n\}$
forms an orthonormal basis with respect to
$\Vert\cdot \Vert_{s,\nu}$. So, to
construct this system one can  orthogonalize the system of functions
$$
h_n(\zeta)=\sum_{l\ge n}c_{l,n}b^{l+1}(\zeta)K^{\alpha\mu^{-(l+1)}}(\zeta).
$$
In other words, $h_n$ is an extremal function of the problem: 
$$
\sup\{\vert \tilde h(0)\vert^2: h=b^{(n+1)}\tilde h,\ \tilde h\in
A^2_1(\Gamma,\alpha\mu^{-(n+1)}),\ 
\Vert  h\Vert^2_{s,\nu}\le1\}.\tag 3.15
$$

Let us compare the extremal problem (3.15) and the extremal problem
$$
\sup\{\vert \tilde h(0)\vert^2: h=b^{(n+1)}\tilde h,\ \tilde h\in
A^2_1(\Gamma,\alpha\mu^{-(n+1)}),\ 
\Vert  h\Vert^2\le1\},
$$
whose extremal function is evidently of the form
$h=b^{n+1}K^{\alpha\mu^{-(n+1)}}$. 
We shall demonstrate in the next sections that
$$
\sum_{\zeta_l\in Z}\left\vert{b^{n+1}(\zeta_l)
K^{\alpha\mu^{-(n+1)}}(\zeta_l)}\right\vert^2{\nu(\zeta_l)}
\to 0,\ n\to\infty,
\tag 3.16
$$
and
$$
P_-(\alpha^{-1})\left\{\bar s b^{n+1} K^{\alpha\mu^{-(n+1)}}\right\}\to 0,\quad
n\to
\infty,\tag 3.17
$$
where $P_-(\alpha^{-1})$ is the orthogonal projection onto
$L^2_{dm\vert\Bbb E}\ominus A^2_1(\Gamma,\alpha^{-1})$. Now (3.17) says that
$\langle \bar s b^{n+1} K^{\alpha\mu^{-(n+1)}},\overline{t b^{n+1}
K^{\alpha\mu^{-(n+1)}}}\rangle\to 0$,
$n\to\infty$. Thus, using also (3.16), we get
$$
\Vert  b^{n+1}
K^{\alpha\mu^{-(n+1)}}\Vert_{s,\nu}\sim \Vert b^{n+1}
K^{\alpha\mu^{-(n+1)}}\Vert.
$$ 
Similarly, $\langle  b^{n+1} K^{\alpha\mu^{-(n+1)}},h_n\rangle_{s,\nu}\sim
\langle  b^{n+1}
K^{\alpha\mu^{-(n+1)}},h_n\rangle$. It implies by (3.14)
$$
p_0\dots p_n\phi(0)\sim K^{\alpha\mu^{-(n+1)}}(0),
$$
and
$$
h_n\sim b^{n+1}K^{\alpha\mu^{-(n+1)}},
$$ 
which gives by (3.3) an asymptotic relation for $P_n$.
Recall that the asymptotic relations have been derived under the assumptions
 (3.11), (3.12). In
the next two sections we will present a wide  class of measures for which 
these assumptions are satisfied.  Approximating a
given measure by a sequence of such measures 
we are able to prove  the Main Theorem.


\head 4. Pure--point spectrum in the gaps
\endhead

In this section we investigate the pure--point spectrum of 
$\sigma$ and $\sigma^{(\tau)}$ on $\Bbb R\setminus E$.

Our goal is to present a sufficient condition such that the
associated point measures
$\nu$ and $\nu^{(\tau)}$
 satisfy  (3.12).
\proclaim{Lemma 4.1} Let $r(z)$ be a function of the form (2.1),
such that $\log\sigma'_{a.c.}\circ z\in L^1$.
Assume that the set of zeros and poles $X\cup X^{(\tau)}$ of $r(z)$ satisfies
$$
\sup_{y_l\in X\cup X^{(\tau)}}\sum_{j\not=l}G(y_j,y_l)\le-\log\delta\quad
(0<\delta<1).\tag 4.1
$$
Then there exists a constant $C(E,\delta)<\infty$, such that
$\forall\beta\in\Gamma^*$
$$
\sum_{x_k\in X}\left\vert\frac g \psi\right\vert^2(x_k){\sigma_k}
\le C(E,\delta)\Vert
g\Vert^2\tag 4.2
$$
and
$$
\sum_{x_k^{(\tau)}\in
X^{(\tau)}}\left\vert
\frac g \phi\right\vert^2(x_k^{(\tau)}){\sigma_k^{(\tau)}}
\le C(E,\delta)\Vert
g\Vert^2,\tag 4.3
$$
where $g\in A_1^2(\Gamma,\beta)$,
$\overline{g(\bar\zeta)}=g(\zeta)$, $g(0)=0$.
\endproclaim
\demo{Proof} We will prove (4.2), (4.3) could be proved in the same way. The
proof is based on the following Theorem (see [6, 15]): if $E$ is a
homogeneous set,
$\bar\Bbb C\setminus E\equiv \Bbb D/\Gamma$, then for any $\beta\in\Gamma^*$
there exists a $w\in H^\infty(\Gamma,\beta)$,
$\Vert w\Vert\le 1$, such that $\vert w(\zeta)\vert\ge C(E)$.
We  may assume, that $\overline{w(\bar\zeta)}=
w(\zeta)$.

Let $\{X^N\}$ be an exhaustion of $X$ by finite sets.
Let $B_N$ be the Blaschke product with
zeros in $\{\zeta_j:\ z(\zeta_j)\in X\cup X^{(\tau)}\setminus X^N\}$,
$B_N(0)>0$. In this case $\overline{B_N(\bar\zeta)}=B_N(\zeta)$. We note, that
at any point
$\zeta_l:\ z(\zeta_l)\in X^N$, we have $\vert B_N(\zeta_l)\vert\ge\delta$.

Let $w$ be a function from the above cited Theorem, such that the character of
the function $B_N w g$ equals $\alpha$. The function $\frac{(B_N w
g)^2}{\phi\psi}z'$  has only a finite number of poles
(all zeros  of the product $\psi\phi$, except a finite number of them, were 
included in $B_N$).
Thus, we are able to apply the (DCT):
$$
\int_{\Bbb E}\frac{(B_N w
g)^2}{\phi\psi}z'\frac{dt}{2\pi i}=
\sum_{z(\zeta_l)\in X^N}\left(\frac{B_N w
g}{\psi}\right)^2\frac\psi{\phi'} z'.
$$
The function $\left(\frac{B_N wg}{\psi}\right)$ is automorphic and 
real--valued on the real axis, therefore it defines a function on
$\bar\Bbb C\setminus E$ which is real--valued on $\Bbb R\setminus E$.
So, we can write square--module instead of square. Then, we get
$$
\sum_{z(\zeta_l)\in X^N}\left(\frac{B_N w
g}{\psi}\right)^2\sigma_l\ge \delta^2 C(E)^2
\sum_{z(\zeta_l)\in X^N}\left\vert\frac{
g}{\psi}\right\vert^2\sigma_l.
$$
On the other hand, taking into account, that 
$\left\vert\frac{ z'}{\phi\psi} \right\vert=
\left\vert\frac{\psi\bar\phi-\phi\bar\psi }{\phi\psi}
\right\vert\le 2$, we get
$$
\left\vert\int_{\Bbb E}\frac{(B_N w
g)^2}{\phi\psi}z'\frac{dt}{2\pi i}\right\vert\le 2
\int_{\Bbb E}\vert
g\vert^2\,dm
$$
Hence the lemma is proved.\qed
\enddemo

\proclaim{Lemma 4.2} Let $\nu$ be a measure  of the form (3.9)
with the property
$$
\left\{\sum_{\zeta_l\in Z}\vert g(\zeta_l)\vert^2\nu(\zeta_l)\right\}^{1/2}
\le C\Vert
g\Vert,\quad g\in A_1^2(\Gamma,\beta).
$$
Then
$$
\left\{\sum_{\zeta_l\in Z}\vert (b^n
K^{\alpha\mu^{-n}})(\zeta_l)\vert^2\nu(\zeta_l)
\right\}^{1/2}\to 0,\quad n\to \infty.
$$
\endproclaim

\demo{Proof} For fixed $\epsilon>0$ let
$$
\Gamma^*=\bigcup_{j=1}^{l(\epsilon)}\{\beta:\ \text{dist}(\beta,\beta_j)\le
\eta(\epsilon)\}
$$
be a finite covering of $\Gamma^*$, such that
$$
2\left\vert1-\frac{\Delta^{\beta_j^{-1}\beta}(0)K^{\beta_j}(0)}
{K^\beta(0)}\right\vert\le\epsilon^2,
\quad
\text{dist}(\beta,\beta_j)\le\eta(\epsilon).
$$
It means, that
$$
\Vert(\Delta^{\beta_j^{-1}\beta}K^{\beta_j}) - K^\beta\Vert^2
\le 1+1-2\frac{\Delta^{\beta_j^{-1}\beta}(0)K^{\beta_j}(0)}
{K^\beta(0)}
\le\epsilon^2,
\quad
\text{dist}(\beta,\beta_j)\le\eta(\epsilon).\tag 4.4
$$
For a fixed $\beta_j$ the sum
$$
\sum_{\zeta_l\in Z}\vert 
K^{\beta_j}(\zeta_l)\vert^2\nu(\zeta_l)<\infty
$$
is a convergent majorant for the series
$$
\sum_{\zeta_l\in Z}\vert (b^n
K^{\beta_j})(\zeta_l)\vert^2\nu(\zeta_l).
$$
Therefore, it tends to 0 as $n$ tends to infinity. It means, that
there exists an $n_0$ such that
$$
\left\{\sum_{\zeta_l\in Z}\vert (b^n
K^{\beta_j})(\zeta_l)\vert^2\nu(\zeta_l)\right\}^{1/2}\le \epsilon,
\ 
\forall n>n_0,\ 1\le j\le l(\epsilon).\tag 4.5
$$

Let, now, $n>n_0$ and
$\beta_j:$ $\text{dist}(\beta_j,\alpha\mu^{-n})\le\eta(\epsilon)$. Then,
due to (4.4) and (4.5) we get
$$\align
&\left\{\sum_{\zeta_l\in Z}\vert (b^n
K^{\alpha\mu^{-n}})(\zeta_l)\vert^2\nu(\zeta_l)
\right\}^{1/2}\\ \le&
\left\{\sum_{\zeta_l\in Z}\vert 
(b^n
K^{\alpha\mu^{-n}})(\zeta_l)
-(b^n
K^{\beta_j}\Delta^{\beta_j^{-1}\alpha\mu^{-n}})(\zeta_l)\vert^2\nu(\zeta_l)
\right\}^{1/2}\\+&
\left\{\sum_{\zeta_l\in Z}\vert (b^n
K^{\beta_j}\Delta^{\beta_j^{-1}\alpha\mu^{-n}})(\zeta_l)\vert^2\nu(\zeta_l)
\right\}^{1/2}\\ \le& C
\Vert K^{\alpha\mu^{-n}}
-
K^{\beta_j}\Delta^{\beta_j^{-1}\alpha\mu^{-n}}\Vert+
\left\{\sum_{\zeta_l\in Z}\vert (b^n
K^{\beta_j})(\zeta_l)\vert^2\nu(\zeta_l)
\right\}^{1/2}\\ \le&(C+1)\epsilon.
\endalign
$$
\qed
\enddemo

Combining these two lemmas we have the following proposition.
\proclaim{Lemma 4.3} Let $r(z)$ be a function of the form (2.1),
such that $\log\sigma'_{a.c.}\circ z\in L^1$.
Assume that the system of zeros and poles $X\cup X^{(\tau)}$ of $r(z)$ satisfies
(4.1).
Then 
$$
\sum_{\zeta_l\in Z}\left\vert{b^n(\zeta_l) K^{\alpha\mu^{-n}}}(\zeta_l)
\right\vert^2{\nu(\zeta_l)}
\to 0,\ n\to\infty,
$$
and
$$
\sum_{\zeta_l^{(\tau)}\in Z^{(\tau)}}\left\vert{b^n(\zeta_l^{(\tau)})
K^{\alpha\mu^{-n}}}(\zeta_l^{(\tau)})
\right\vert^2{\nu^{(\tau)}(\zeta_l^{(\tau)})}
\to 0,\ n\to\infty,
$$
where $\nu$ and $\nu^{(\tau)}$ are defined by (3.7) and (3.8).
\endproclaim


\head 5. Enclosure in $A_1^2(\Gamma, \alpha)$. Absolutely continuous
spectrum\endhead

In this section we present  sufficient conditions for
the map $P\to h(\zeta, P;\sigma)$ to be a bounded map
from $L^2_{d\sigma}$ to $A_1^2(\Gamma,\alpha)$. 

With a given set of interlacing points $X\cup X^{(\tau)}$
from $\Bbb R\setminus E$,
$E=[b_0,a_0]\setminus\cup_{j\ge 1}(a_j,b_j)$, we associate
a special function $r_0(z)=r_0(z;X\cup X^{(\tau)})$ of the form (2.8)  with the given
set of  poles and zeros  by putting $u(x)=\pi/2$. Let $d\sigma_0(x)=
d\sigma_0(x;X\cup X^{(\tau)})$ be the measure associated with this function,
$$
1+\int\frac{d\sigma_0(x;X\cup X^{(\tau)})}{x-z}=r_0(z;X\cup X^{(\tau)})=
\prod_l\frac{z-x_l^{(\tau)}}{z-x_l}\sqrt{\prod_{j\ge 0}
\frac{z-a_j}{z-b_j}}.\tag 5.1
$$

\proclaim{Lemma 5.1} Let $X\cup X^{(\tau)}\subset\Bbb R\setminus E$ be a  given
set of interlacing points, such that $X$ satisfies condition (2.5). Let the
measure $d\sigma_0(z)=d\sigma_0(z;X\cup X^{(\tau)})$ be defined by (5.1).

Then
$P\to h(\zeta, P;\sigma_0)$ is a bounded map.
Moreover
$$
\int_{\Bbb E}\vert h\vert^2 \,dm\le 2
\int P^2\,d\sigma_0.
$$
\endproclaim
\demo{Proof} First note, that $d\sigma_0$ is absolutely continuous
on $E$. In fact, $r_0(x+i0)$ takes pure imaginary values on $E$,
so $\vert r_0(x)\vert dx=\text{Im}\ r_0(x) dx$ and we can apply
Lemma 2.4.  Using Lemma 2.2, we present $r_0(z;X\cup X^{(\tau)})$ in the form
$r_0=\psi_0/\phi_0$ with the property (2.7).
We put $s_0=\phi_0/(\overline{t\phi_0})$,
$s_0^{(\tau)}=\psi_0/(\overline{t\psi_0})$.

As a consequence of (3.5), (3.6), we obtain the estimate (see also Lemma 2.1):
$$
\align
&\langle
\bmatrix 1& -s_0\\ 1& -s_0^{(\tau)}\endbmatrix\bmatrix h(t)\\\bar t h(\bar
t)\endbmatrix,
\bmatrix 1& -s_0\\ 1& -s_0^{(\tau)}\endbmatrix\bmatrix h(t)\\\bar t h(\bar
t)\endbmatrix\rangle_{L^2_{dm\vert\Bbb E}}\\
=& \int_{\Bbb E} P^2\left\vert\frac{t z'}{\phi_0}\right\vert^2\,dm
+\int_{\Bbb E} (P^{(\tau)})^2\left\vert\frac{t z'}{\psi_0}\right\vert^2\,dm\\
\le& 2\int_E P^2\,d\sigma_{0}+2\int_E  (P^{(\tau)})^2\,d\sigma_{0}^{(\tau)}\\
=& 4\int P^2\,d\sigma_0.\tag 5.2
\endalign
$$
Since
$$
\bmatrix 1& -s_0\\ 1& -s_0^{(\tau)}\endbmatrix^*\bmatrix 1& -s_0\\ 1&
-s_0^{(\tau)}\endbmatrix =2\bmatrix 1& -\frac{s_0+s_0^{(\tau)}} 2\\
-\frac{\overline{s_0+s_0^{(\tau)}}} 2& 1\endbmatrix,
$$
and $\left\vert\frac{s_0+s_0^{(\tau)}}
2\right\vert^2+\left\vert\frac{s_0-s_0^{(\tau)}} 2\right\vert^2=1$, we have
$$
2\bmatrix \left\vert\frac{s_0+s_0^{(\tau)}} 2\right\vert^2& -\frac{s_0+s_0^{(\tau)}}
2\\ -\frac{\overline{s_0+s_0^{(\tau)}}} 2& 1\endbmatrix+
2\bmatrix \left\vert\frac{s_0-s_0^{(\tau)}} 2\right\vert^2& 0\\
0& 0\endbmatrix\ge
2\left\vert\frac{s_0-s_0^{(\tau)}} 2\right\vert^2
\bmatrix 1& 0\\
0& 0\endbmatrix.\tag 5.3
$$
At last, 
$$
\left\vert\frac{s_0-s_0^{(\tau)}} 2\right\vert^2=
\left\vert\frac{1-r_0/\bar r_0} 2\right\vert^2=\sin^2(\arg
r_0)=1.
\tag 5.4
$$

Combining (5.2), (5.3) and (5.4), we get
$$
2
\int_{\Bbb E}\vert h\vert^2 \,dm\le 4
\int P^2\,d\sigma_0,
$$
which proves the assertion.\qed
\enddemo
\proclaim{Remark} Let us note, that in the same way it can be proved that
the map $P\to h(\zeta, P,\sigma)$ is bounded if the associated
Stieltjes function $r(z)$ possesses the property
$$
\delta\le\arg r(x+i0)\le\pi-\delta, \quad x\in E,
$$
with some $\delta>0$. In this case
$$
\int_{\Bbb E}\vert h\vert^2 \,dm\le \frac 2{\sin^2\delta}
\int P^2\,d\sigma.
$$
\endproclaim

\proclaim{Lemma 5.2} Let $X\cup X^{(\tau)}$ be a  given set of interlacing
points, such that $X$ satisfies condition (2.5). Let the
measure $d\sigma_0(z)=d\sigma_0(z;X\cup X^{(\tau)})$ be defined by (5.1).
Assume that the measure $d\sigma$ is equivalent to $d\sigma_0$, i.e.
$$
C_1 d\sigma_0\ge d\sigma\ge C_2 d\sigma_0.
$$
Then
$P\to h(\zeta, P;\sigma)$ is a bounded map.
Moreover
$$
\int_{\Bbb E}\vert h\vert^2 \,dm\le 2 \frac{C_1}{C_2}
\int P^2\,d\sigma.
$$
\endproclaim
\demo{Proof} 
Let us represent $d\sigma$ in the form
$$
d\sigma=\frac 1{w^2(x)}\,d\sigma_0,
$$
where $1/C_1 \le w^2(x)\le 1/C_2. 
$ Define $ W(\zeta)$ as an outer function with the given
modulus of boundary values, $\vert W(t)\vert^2=w^2(z(t))$,
$t\in\Bbb T$. In this case the functions $\phi$ and $\phi_0$ associated
with $\sigma$ and $\sigma_0$, respectively, are related in the following
way $\phi=W\phi_0$, and the function $h(\zeta,P;\sigma) $ can be represented
in the form
$$
h(\zeta,P;\sigma)=h(\zeta,\frac{P}{w^2};\sigma_0) W(\zeta).
$$
Therefore,
$$
\Vert h(\zeta,P;\sigma)\Vert^2\le\frac 1{C_2}
\left\Vert h(\zeta,\frac{P}{w^2};\sigma_0)\right\Vert^2\le\frac 2{C_2}
\left\Vert\frac{P}{w^2}\right\Vert^2_{L^2_{d\sigma_0}}\le 2\frac {C_1}{C_2}
\left\Vert P\right\Vert^2_{L^2_{d\sigma}}.$$
\qed
\enddemo
\bigskip

Let us prove (3.17)  not only for an individual function $f=s$,
but for a family of functions $f_\xi$ depending continuously on a
parameter $\xi$.

\proclaim{Lemma 5.3}
Let $\xi\mapsto f_\xi(t)$
be a
$L^\infty(\alpha^2)$--valued continuous function on a compact set
$\Xi$, 
$\Vert f_\xi\Vert_{ L^\infty(\alpha^2)}\le 1$, $\xi\in\Xi$.
Then 
$$
P_-(\alpha^{-1})\left\{\bar f_\xi b^n K^{\alpha\mu^{-n}}\right\}\to 0,\quad
n\to
\infty,\ \text{uniformly on}\  \xi\in\Xi,
$$
where $P_-(\alpha^{-1})$ is the orthogonal projection onto
$L^2_{dm\vert\Bbb E}\ominus A^2_1(\Gamma,\alpha^{-1})$.
 \endproclaim

\demo{Proof} For fixed $\epsilon>0$ let
$$
\Gamma^*=\bigcup_{j=1}^{l(\epsilon)}\{\beta:\ \text{dist}(\beta,\beta_j)\le
\eta(\epsilon)\}
$$
be the same  finite covering of $\Gamma^*$ as in the proof
of Lemma 4.2, i.e.:
$$
\Vert(\Delta^{\beta_j^{-1}\beta}K^{\beta_j}) - K^\beta\Vert
\le\epsilon,
\quad
\text{dist}(\beta,\beta_j)\le\eta(\epsilon).
$$

Let
$$
\Xi=\bigcup_{j'=1}^{l'(\epsilon)}\{\xi:\ \text{dist}(\xi,\xi_{j'})\le
\eta'(\epsilon)\}
$$
be a finite covering of $\Xi$, such that
$$
\Vert
f_\xi-f_{\xi_{j'}}
\Vert\le\epsilon,
\quad
\text{dist}(\xi,\xi_{j'})\le\eta'(\epsilon).
$$

For fixed $\beta$ and $\xi$ one can find $n_0$ such that
$$
\Vert P_-^{n}(\alpha^{-2}\beta)\bar f_\xi K^\beta\Vert\le \epsilon,\ 
\forall n>n_0,
$$
where $P_-^n(\alpha)$ is the orthogonal projection onto
$L^2_{dm\vert\Bbb E}\ominus b^{-n}A^2_1(\Gamma,\alpha\mu^n)$.
Therefore, there exists $n_0$ such that
$$
\Vert P_-^{n}(\alpha^{-2}\beta_j)\overline{ f_{\xi_{j'}}} K^{\beta_j}\Vert\le
\epsilon,\ 
\forall n>n_0,\ 1\le j\le l(\epsilon),\ 1\le j'\le l'(\epsilon).
$$

From now on let $n>n_0=n_0(\epsilon)$. Let $\beta_j:$ 
$\text{dist}(\beta_j,\alpha\mu^{-n})\le\eta(\epsilon)$
and $\xi_{j'}:$ 
$\text{dist}(\xi_{j'},\xi)\le\eta'(\epsilon)$. 
For $h\in L^2_{dm\vert\Bbb E}\ominus A^2_1(\Gamma,\alpha^{-1})$, we write
$$
\align
&\langle b^n K^{\alpha\mu^{-n}}, f_\xi h\rangle=
\langle b^n K^{\alpha\mu^{-n}}, (f_\xi-f_{\xi_{j'}})h\rangle\\
+&\langle b^n (K^{\alpha\mu^{-n}}-\Delta^{\alpha\mu^{-n}\beta_j^{-1}}
K^{\beta_j}), f_{\xi_{j'}} h\rangle +
\langle b^n \Delta^{\alpha\mu^{-n}\beta_j^{-1}}
K^{\beta_j}, f_{\xi_{j'}} h\rangle.
\endalign
$$
Then
$$
\vert\langle b^n K^{\alpha\mu^{-n}},
(f_\xi-f_{\xi_{j'}}) h\rangle\vert
\le \Vert (f_\xi-f_{\xi_{j'}})\Vert \Vert h\Vert
\Vert K^{\alpha\mu^{-n}}\Vert\le \epsilon \Vert h\Vert,
$$
and
$$
\vert\langle b^n (K^{\alpha\mu^{-n}}-\Delta^{\alpha\mu^{-n}\beta_j^{-1}}
K^{\beta_j}), f_{\xi_{j'}} h\rangle\vert
\le \Vert f_{\xi_{j'}}\Vert \Vert h\Vert
\Vert K^{\alpha\mu^{-n}}-\Delta^{\alpha\mu^{-n}\beta_j^{-1}}
K^{\beta_j}\Vert\le \epsilon\Vert h\Vert.
$$
And for the last term we have
$$
\align
&\vert\langle b^n \Delta^{\alpha\mu^{-n}\beta_j^{-1}}
K^{\beta_j}, f_{\xi_{j'}} h\rangle\vert
=\vert\langle \overline{f_{\xi_{j'}}} K^{\beta_j}, b^{-n}
\overline{\Delta^{\alpha\mu^{-n}\beta_j^{-1}}}
 h\rangle\vert\\
=&\vert\langle P_-^{n}(\alpha^{-2}\beta_j) \overline{ f_{\xi_{j'}}} K^{\beta_j},
b^{-n}
\overline{\Delta^{\alpha\mu^{-n}\beta_j^{-1}}}
 h\rangle\vert
\le\Vert  h\Vert 
\Vert P_-^{n}(\alpha^{-2}\beta_j)  \overline{f_{\xi_{j'}}} K^{\beta_j}\Vert\le
\epsilon\Vert h\Vert.
\endalign
$$
Therefore,
$
\vert\langle  P_-(\alpha^{-1})\left\{ \bar f_\xi b^n
K^{\alpha\mu^{-n}}\right\},  h\rangle\vert\le 3\epsilon\Vert h\Vert.
$ Putting $h=P_-(\alpha^{-1})\left\{ \bar f_\xi b^n
K^{\alpha\mu^{-n}}\right\}$, we get
$\Vert P_-(\alpha^{-1})\left\{ \bar f_\xi b^n
K^{\alpha\mu^{-n}}\right\}\Vert\le 3\epsilon.$\qed
\enddemo


\head{ 6. Main Theorem}\endhead

\proclaim{Theorem} Let $E$ be a homogeneous set
and $X\subset\Bbb R\setminus E$ be
a set of points which can accumulate only to the set $E$. Let $\sigma$ be 
a positive measure with the support $E\cup X$. Assume that $\log\sigma'_{a.c.}(z(t))\in
L^1$ and that the set of poles $X$ and zeros $X^{(\tau)}$ of the Stieltjes
function
$$
r(z)=1+\int\frac{d\sigma}{x-z},
$$
satisfies the condition
$$
\sup_{y_l\in X\cup X^{(\tau)}}\sum_{j\not=l} G(y_j,y_l)<\infty.\tag 6.1
$$
Then the minimum deviation and the orthonormal polynomials
$P_n(z,\sigma)=\frac{z^n}{p_0\dots p_n}+\dots$ have the
following asymptotic behavior ($n\to\infty$)
$$
\frac{K^{\alpha\mu^{-(n+1)}}(0)}{p_0\dots p_n\phi(0)}\to 1,\tag 6.2
$$
$$
\overline{\left(\frac{-tz'}{\phi}\right)}P_n-\{b^{n+1}K^{\alpha\mu^{-(n+1)}}-s
\overline{(tb^{n+1}K^{\alpha\mu^{-(n+1)}})}\}\to 0\quad
\text{in}\  L^2_{dm\vert\Bbb E},\tag 6.3
$$
and
$$
\left\vert P_n b^n-\frac{\phi}{(-z')b^2}K^{\alpha^{-1}\mu^{n+2}}\right\vert
\to 0\tag 6.4
$$
uniformly on each compact subset of $\bar\Bbb C\setminus E$.
Here $\phi$ is chosen  and $\alpha$ is given as in Lemma 2.2 and
$s=\phi/\overline{(t\phi)}$.
\endproclaim

The proof of the theorem will be divided into several steps. The main part deals
with the statement (6.2). First we show an upper estimate.
\proclaim{Lemma 6.1} Under the assumptions of the previous Theorem, we have
$$
\overline\lim_{n\to\infty}\frac{K^{\alpha\mu^{-(n+1)}}(0)}{p_0\dots
p_n\phi(0)}\le 1.\tag 6.5
$$
Furthermore, (6.2) implies (6.3).
\endproclaim

\demo{Proof} Consider the norm of the function in (6.3),
$$
\aligned
&\left\Vert\overline{\left(\frac{-tz'}{\phi}\right)}P_n-\{b^{n+1}K^{\alpha\mu^{-(n+1)}}-s
\overline{(tb^{n+1}K^{\alpha\mu^{-(n+1)}})}\}\right\Vert^2\\
=&\int_{\Bbb E}P_n^2\left\vert\frac{tz'}{\phi}\right\vert^2\,dm-
2\text{Re}\ \left\langle
b^{n+1}K^{\alpha\mu^{-(n+1)}}-s
\overline{(tb^{n+1}K^{\alpha\mu^{-(n+1)}})},\overline{\left(\frac{tz'}{\phi}\right)}P_n
\right\rangle\\+&
\Vert b^{n+1}K^{\alpha\mu^{-(n+1)}}-s
\overline{(tb^{n+1}K^{\alpha\mu^{-(n+1)}})}\Vert^2.
\endaligned\tag 6.6
$$
To prove (6.5)
we only use the fact that the norm is non-negative. From the estimate we get it
follows immediately, that (6.2) implies (6.3).

For the first term in (6.6) we have an estimate
$$
\int_{\Bbb E}P_n^2\left\vert\frac{tz'}{\phi}\right\vert^2\,dm=
2\int P_n^2\,d\sigma_{a.c.}\le2\int P_n^2\,d\sigma=2.\tag 6.7
$$
Due to Lemma 5.3,
$$
\langle
b^{n+1}K^{\alpha\mu^{-(n+1)}},s
\overline{(tb^{n+1}K^{\alpha\mu^{-(n+1)}})}\rangle\to 0,\ n\to\infty,
$$
so,
$$
\align
\Vert b^{n+1}K^{\alpha\mu^{-(n+1)}}-s
\overline{(tb^{n+1}K^{\alpha\mu^{-(n+1)}})}\Vert^2=& 1+1-2\text{Re}
\langle
b^{n+1}K^{\alpha\mu^{-(n+1)}},s
\overline{(tb^{n+1}K^{\alpha\mu^{-(n+1)}})}\rangle\\ \to& 2,\quad n\to\infty.
\tag 6.8
\endalign
$$
Note, also, that $\overline{\left(\frac{-tz'}{\phi}\right)}P_n$ is orthogonal
to any function of the form $g(t)+s(t)\bar t g(\bar t)$. Therefore,
$$
\left\langle
b^{n+1}K^{\alpha\mu^{-(n+1)}}-s
\overline{(tb^{n+1}K^{\alpha\mu^{-(n+1)}})},\overline{\left(\frac{-tz'}{\phi}\right)}P_n
\right\rangle= 2
\left\langle
b^{n+1}K^{\alpha\mu^{-(n+1)}},\overline{\left(\frac{-tz'}{\phi}\right)}P_n
\right\rangle.
$$
To evaluate this scalar product we apply the (DCT),
$$\aligned
\left\langle
b^{n+1}K^{\alpha\mu^{-(n+1)}},\overline{\left(\frac{-tz'}{\phi}\right)}P_n
\right\rangle=&
\int_{\Bbb E}
b^{n+1}K^{\alpha\mu^{-(n+1)}}\left(\frac{-tz'}{\phi}\right)P_n
\frac{dt}{2\pi i t}\\
=&
\int_{\Bbb E}
K^{\alpha\mu^{-(n+1)}}\left(\frac{-b^2z'}{\phi}\right)(b^n P_n)
\frac{dt}{2\pi i b}\\
=&
\frac 1{b'(0)}\left\{
K^{\alpha\mu^{-(n+1)}}\left(\frac{-b^2z'}{\phi}\right)(b^n P_n)
\right\}(0)\\
+&
\sum_{\zeta_l\in Z}
\left\{b^{n+1}K^{\alpha\mu^{-(n+1)}}\left(\frac{-z'}{\phi'}\right)P_n
\right\}(\zeta_l).
\endaligned\tag 6.9
$$
Since $(bz)(0)=1$, we have
$$
\frac 1{b'(0)}\left\{
K^{\alpha\mu^{-(n+1)}}\left(\frac{-b^2z'}{\phi}\right)(b^n P_n)
\right\}(0)=\frac{
K^{\alpha\mu^{-(n+1)}}(0)}{p_0\dots p_n\phi(0)}.
$$
Let us show that the last term in (6.9) tends to 0.
Indeed,
$$
\align &
\sum_{\zeta_l\in
Z}\left\vert b^{n+1}K^{\alpha\mu^{-(n+1)}} P_n\frac{(-z')}{\phi'}\right\vert=
\sum_{\zeta_l\in Z}\vert
P_n\vert\left\vert\frac{b^{n+1}K^{\alpha\mu^{-(n+1)}}}{\psi}\right\vert 
\frac{(- z')\psi}{\phi'}\\
\le&\sqrt{\sum_{\zeta_l\in Z}\vert P_n\vert^2\sigma_l}
\sqrt{\sum_{\zeta_l\in
Z}\left\vert\frac{b^{n+1}K^{\alpha\mu^{-(n+1)}}}{\psi}\right\vert^2 
\sigma_l}.
\endalign
$$
Since
$$
{\sum_{x_l\in X}\vert P_n\vert^2\sigma_l}\le
{\int\vert
P_n\vert^2\,d\sigma}=1,
$$
using Lemma 4.3 and the definition of the measure $\nu$, we get
$$
\sum_{\zeta_l\in Z}
\left\{b^{n+1}K^{\alpha\mu^{-(n+1)}}\left(\frac{-z'}{\phi'}\right)P_n
\right\}(\zeta_l)\to 0.\tag 6.10
$$
Substituting (6.7)...(6.10) in (6.6) we obtain
$$\align
&2-4\frac{
K^{\alpha\mu^{-(n+1)}}(0)}{p_0\dots p_n\phi(0)}+2+o(1)\\ \ge &
\left\Vert\overline{\left(\frac{-tz'}{\phi}\right)}P_n-\{b^{n+1}K^{\alpha\mu^{-(n+1)}}-s
\overline{(tb^{n+1}K^{\alpha\mu^{-(n+1)}})}\}\right\Vert^2\ge 0,
\endalign
$$
and thus the lemma is proved.\qed
\enddemo
\proclaim{Lemma 6.2}Assume  that
$$
\Vert h(\zeta,P;\sigma)\Vert\le C\Vert P\Vert_{L^2_{d\sigma}},
$$ 
and that the assumptions of the previous Theorem are satisfied.
Then
$$
\lim_{n\to\infty}\frac{K^{\alpha\mu^{-(n+1)}}(0)}{p_0\dots
p_n\phi(0)}= 1.
$$
\endproclaim
\demo{Proof}
Let $h_n(\zeta)=h_n(\zeta,P_n;\sigma)$. Then
$$
\overline{\left(\frac{-tz'}{\phi}\right)}P_n(z(t))=h_n-s\bar th_n(\bar t),
$$
where $s=\phi/\overline{(t\phi)}$. Multiplying this identity
by $b^{n+1}K^{\alpha\mu^{-(n+1)}}$, we get (see (6.9), (6.10))
$$
\frac{K^{\alpha\mu^{-(n+1)}}(0)}{p_0\dots
p_n\phi(0)}+o(1)=
\langle h_n,b^{n+1}K^{\alpha\mu^{-(n+1)}}\rangle
-\langle s\bar th_n(\bar
t),b^{n+1}K^{\alpha\mu^{-(n+1)}}\rangle
\tag 6.11
$$
Using the reproducing property of $K^{\alpha\mu^{-(n+1)}}$ and (3.14) we get
$$
\langle h_n,b^{n+1}K^{\alpha\mu^{-(n+1)}}\rangle=
\frac
{p_0\dots
p_n\phi(0)}
{K^{\alpha\mu^{-(n+1)}}(0)},
$$
and due to Lemma 5.3 we have
$$
\aligned
\vert\langle s\bar t h_n(\bar
t),b^{n+1}K^{\alpha\mu^{-(n+1)}}\rangle\vert=&
\vert\langle  \bar t h_n(\bar t),P_-(\alpha^{-1})
\{\bar s b^{n+1}K^{\alpha\mu^{-(n+1)}}\}\rangle\vert\\
\le&\Vert h_n\Vert
\Vert P_-(\alpha^{-1})
\{\bar s b^{n+1}K^{\alpha\mu^{-(n+1)}}\}\Vert\to 0,\ n\to\infty.
\endaligned
$$
Thus (6.11) is of the form,
$$
\frac{K^{\alpha\mu^{-(n+1)}}(0)}{p_0\dots
p_n\phi(0)}+o(1)=
\frac{p_0\dots
p_n\phi(0)}
{K^{\alpha\mu^{-(n+1)}}(0)}+o(1)
$$
and the lemma is proved.\qed
\enddemo

\demo{Proof of (6.2) and (6.3)} Let $r_0(z)=r_0(z;X\cup X^{(\tau)})$ be a
function of the form (5.1), associated with the given set of zeros and
poles. For $0<\eta<1$, put
$$
r_\eta(z)=\eta r(z)+(1-\eta)r_0(z).
$$
Note that for an arbitrary $\eta$ it is also the Stieltjes
function with the same set of zeros and poles as $r$.
In what follows all functions and coefficients, related to the
function $r_\eta(z)$ have the same subscript $\eta$. We
need some facts concerning these objects.

If $P_{n,\eta}(z)=\frac{z^n}{p_{0,\eta}\dots p_{n,\eta}}+\dots$
is the orthonormal polynomial with respect to $\sigma_\eta$,
then
$$
\aligned
p^2_{0}\dots p^2_{n}=&\inf_{\{P(z)=z^n+\dots\}}\int
P^2\,d\sigma\\ \le&
\int\vert{p_{0,\eta}\dots p_{n,\eta}}P_{n,\eta}\vert^2\,d\sigma\\ \le&
\frac 1\eta\int\vert{p_{0,\eta}\dots
p_{n,\eta}}P_{n,\eta}\vert^2\,d\sigma_\eta\\ =&
\frac 1\eta{p^2_{0,\eta}\dots p^2_{n,\eta}}.
\endaligned\tag 6.12
$$

Further, since
$$
\frac 1{\vert\phi_{\eta}\vert^2}=
\frac {\eta}{\vert\phi\vert^2}+
\frac {1-\eta}{\vert\phi_{0}\vert^2},
$$
the function $\log{\vert\phi_{\eta}\vert^2}$ converges to
$\log\vert\phi\vert^2 $ in $L^1$, as $\eta\to 1$. Therefore,
$$
\alpha_\eta\to\alpha\quad\text{and}\quad
\phi_\eta(0)\to\phi(0),\quad\eta\to 1.\tag 6.13
$$

The measure $d\sigma^{(\tau)}_\eta$ is absolutely continuous on $E$.
Indeed,
$$
\left\vert \frac1 {r_\eta}\right\vert\le
\left\vert \frac1 {\text{Im}\ r_\eta}\right\vert\le\frac 1{1-\eta}
\left\vert \frac1 {\text{Im}\ r_0}\right\vert\le\frac 1{1-\eta}
\left\vert \text{Im}\ \frac{1} {r_0}\right\vert.
$$
And since $-\text{Im}\ \frac{1} {r_0}\,dx$ is integrable
we can apply Lemma 2.4. Simultaneously we have an estimate
$$
d\sigma_\eta^{(\tau)}\le\frac 1{1-\eta}d\sigma^{(\tau)}_0\quad\text{on}\ E.
$$
Since $\sigma^{(\tau)}_l=1/r'(x_l^{(\tau)})$, and
$$
r'_\eta(x_l^{(\tau)})=\eta r'(x_l^{(\tau)})+(1-\eta)r_0'(x_l^{(\tau)}),
$$
we also have
$$
\sigma_{l,\eta}^{(\tau)}\le\frac 1{1-\eta}\sigma^{(\tau)}_{l,0}.
$$
So,
$$
d\sigma_\eta^{(\tau)}\le\frac 1{1-\eta}d\sigma^{(\tau)}_0.\tag 6.14
$$

Next, let us put
$$
-\frac 1{r_{\eta,\eta_1}}=
-\frac {\eta_1}{r_{\eta}}-\frac {1-\eta_1}{r_0}.
$$
Then, due to (6.14),
$$
(1-\eta_1)d\sigma_0^{(\tau)}
\le d\sigma_{\eta,\eta_1}^{(\tau)}\le
\left[\frac{\eta_1}{1-\eta}+(1-\eta_1)\right]{d\sigma_0^{(\tau)}}.
$$
Now to this measure we can apply Lemma 5.2 and
Lemma 6.2. Therefore,
$$
\lim_{n\to\infty}\frac{K^{\alpha_{\eta,\eta_1}\mu^{-(n+1)}}(0)}
{p_{0,\eta,\eta_1}\dots p_{n,\eta,\eta_1}\psi_{\eta,\eta_1}(0)}=
1.\tag 6.15
$$
We remind that the functions $\psi$ and $\phi$ have the same character
and that $\phi(0)=\psi(0)$. As above in (6.13), since
$\log\vert\psi_{\eta,\eta_1}\vert\to\log\vert\psi_\eta\vert$ in
$L^1$ ($\eta_1\to 1$),
we have
$$
\alpha_{\eta,\eta_1}\to\alpha_\eta\quad\text{and}\quad
\psi_{\eta,\eta_1}(0)\to\psi_\eta(0),\quad\eta_1\to 1.\tag 6.16
$$
At last, as in (6.12) we obtain
$$
p^2_{0,\eta,\eta_1}\dots p^2_{n,\eta,\eta_1}
\ge\eta_1\ {p^2_{0,\eta}\dots p^2_{n,\eta}}
\ge\eta_1\eta\ {p^2_{0}\dots p^2_{n}}.
\tag 6.17
$$

Since $K^\alpha(0)$ depends continuously on $\alpha$,
using (6.13) and (6.16), for arbitrary $\epsilon>0$, we can chose $\eta, \eta_1$
so close to 1, that $\forall n$,
$$
\frac{K^{\alpha\mu^{-(n+1)}}(0)}
{\phi(0)}\ge
\frac{K^{\alpha_{\eta,\eta_1}\mu^{-(n+1)}}(0)}
{\psi_{\eta,\eta_1}(0)}(1-\epsilon).
$$
Using (6.17), we get
$$
\frac{K^{\alpha\mu^{-(n+1)}}(0)}
{p_0\dots p_n\phi(0)}\ge
\sqrt{\eta\eta_1}\frac{K^{\alpha_{\eta,\eta_1}\mu^{-(n+1)}}(0)}
{p_{0,\eta,\eta_1}\dots p_{n,\eta,\eta_1}\psi_{\eta,\eta_1}(0)}(1-\epsilon).
$$
Then, due to (6.15), we get
$$
\underline\lim\Sb{n\to\infty}\endSb\frac{K^{\alpha\mu^{-(n+1)}}(0)}
{p_0\dots p_n\phi(0)}\ge 1.
$$
Together with (6.5) this finishes the proof of (6.2),
and due to Lemma 6.1, (6.3) is also proved.\qed

\enddemo

\demo{Proof of (6.4)}
We rewrite (6.3) in terms of analytic functions.
Taking into account, that 
$\overline{(t K^{\alpha\mu^{-(n+1)}})}=\frac{K^{\alpha^{-1}\mu^{n+2}}} b$,
we have
$$
g_n:=\frac{(-z')}\phi P_n-\left\{b^{-(n+2)}K^{\alpha^{-1}\mu^{n+2}}
-\bar s b^{(n+1)}K^{\alpha\mu^{-(n+1)}}
\right\}\to 0\quad\text{in}\ L^2_{dm\vert\Bbb E}.
$$
In a standard way, using reproducing kernel, we get
$$
\left\vert\left\langle g_n,
\overline{(b^{n+2} B_{X})}
k^{\alpha_{\text{out}}^{-1}\mu^{n+2}}(t,\zeta)\right\rangle
\right\vert
\le\left\Vert
g_n\right\Vert
\sqrt{k^{\alpha_{\text{out}}^{-1}\mu^{n+2}}(\zeta,\zeta)},
\tag 6.18
$$
moreover,
$$
\aligned
&\left\langle
\frac{(-z')}\phi P_n-b^{-(n+2)}K^{\alpha^{-1}\mu^{n+2}}
, \overline{(b^{n+2} B_{X})}
k^{\alpha_{\text{out}}^{-1}\mu^{n+2}}(t,\zeta)\right\rangle
\\
=&\left(
\frac{(- z')b^2}{\phi^{\text{out}}} b^n P_n-
B_{X}K^{\alpha^{-1}\mu^{n+2}}\right) (\zeta).
\endaligned\tag 6.19
$$
Therefore, our goal is to evaluate the scalar product
$$\left\langle
 b^{n+1}K^{\alpha\mu^{-(n+1)}}
, s\overline{(b^{n+2} B_{X})}
k^{\alpha_{\text{out}}^{-1}\mu^{n+2}}(t,\zeta)\right\rangle.
$$
We claim that this product tends to $0$ uniformly with respect to $\zeta$.

Let $w$ be an inner function from $H^\infty(\Gamma)$, for example,
$w=b\Delta^{\mu^{-1}}$. Then
$$
\overline{\left(\frac{w-w(\zeta)}{1-w\overline{w(\zeta)}}
\right)}k^{\alpha_{\text{out}}^{-1}\mu^{n+2}}(t,\zeta)\in L^2_{dm\vert\Bbb E}\ominus
A_1^2(\Gamma,{\alpha_{\text{out}}^{-1}\mu^{n+2}}),\tag 6.20
$$
and
$$
\left\Vert \overline{B_{X_-}b^{n+2}}
\overline{\left(\frac{w-w(\zeta)}{1-w\overline{w(\zeta)}}
\right)}k^{\alpha_{\text{out}}^{-1}\mu^{n+2}}(t,\zeta)
\right\Vert=\sqrt{k^{\alpha_{\text{out}}^{-1}\mu^{n+2}}(\zeta,\zeta)}.\tag 6.21
$$
For a fixed compact $\Xi\subset\bar\Bbb C\setminus E$, put
$f(t,z)=s\frac{w-w(\zeta)}{1-w\overline{w(\zeta)}}$, $z(\zeta)=z\in\Xi$.
This is a $L^\infty(\alpha^2)$--valued continuous function on $\Xi$.
Since
$$\align
&\left\langle
 b^{n+1}K^{\alpha\mu^{-(n+1)}}
, s\overline{(b^{n+2} B_{X})}
k^{\alpha_{\text{out}}^{-1}\mu^{n+2}}(t,\zeta)\right\rangle\\
=&\left\langle
 b^{n+1}K^{\alpha\mu^{-(n+1)}}
,f(t,z)\overline{B_{X}b^{n+2}}
\overline{\left(\frac{w-w(\zeta)}{1-w\overline{w(\zeta)}}
\right)}k^{\alpha_{\text{out}}^{-1}\mu^{n+2}}(t,\zeta)\right\rangle\\
=&\left\langle P_-(\alpha^{-1})\left\{\overline{f(t,z)}
 b^{n+1}K^{\alpha\mu^{-(n+1)}}\right\}
,\overline{B_{X}b^{n+2}}
\overline{\left(\frac{w-w(\zeta)}{1-w\overline{w(\zeta)}}
\right)}k^{\alpha_{\text{out}}^{-1}\mu^{n+2}}(t,\zeta)\right\rangle,
\endalign
$$
(see (6.20)) we can apply  Lemma 5.3 . Together with (6.18), (6.19), (6.21),
we obtain the estimate: $\forall\epsilon$ $\exists n_0$, such that
$$
\left\vert\left(
\frac{(-z')b^2}{\phi^{\text{out}}} b^n P_n\right)(\zeta)-
\left(B_{X}K^{\alpha^{-1}\mu^{n+2}}\right) (\zeta)\right\vert\le\epsilon
\sup_{\alpha\in\Gamma^*}\sqrt{k^{\alpha}(\zeta,\zeta)}
,\quad n\ge n_0,\ z(\zeta)\in\Xi.
$$
Note, that 
$$
\sup_{\alpha\in\Gamma^*}\frac{\sqrt{k^{\alpha}(\zeta,\zeta)}
\vert\phi^{\text{out}}(\zeta)\vert}{\vert (b^2 z')(\zeta)\vert}
$$ 
 defines a function in $\bar\Bbb C\setminus
E\equiv\Bbb D/\Gamma$, which is uniformly bounded on $\Xi$. Therefore,
$$
\left\vert
 b^n P_n -
\frac{B_{X}\phi^{\text{out}}}{(-z')b^2}K^{\alpha^{-1}\mu^{n+2}}
\right\vert\le C(\Xi)\epsilon
$$
and the theorem is proved.\qed
\enddemo

\demo{Proof of Corollary 0.1}
In this case $X=\emptyset$, and $X^{(\tau)}$ contains at most one point in
each gap. Therefore, due to the lemma of Jones and Marshall (see Sect. 1) condition (6.1)
holds automatically. Moreover $\phi$ is an outer function. So, substituting
(2.9) in (6.2), (6.3) and (6.4) we get (0.18), (0.20) and (0.19)
respectively.\qed
\enddemo

\proclaim{Corollary 6.1} Under the assumptions of the previous Theorem, the
recurrence coefficients $\{p_n\}$, $\{q_n\}$ have the following
asymptotic behavior
$$
p_n-\Cal P(\alpha^{-1}\mu^{n+1})\to 0,\quad n\to\infty,\tag 6.22
$$
and
$$
q_{n-1}-\Cal Q(\alpha^{-1}\mu^{n+1})\to 0,\quad n\to\infty.\tag 6.23
$$
Moreover, $\{p_n\}$ and $\{q_n\}$ are limit (uniformly) almost periodic.
\endproclaim
\demo{Proof} As it was already shown, (we still assume, that
$(bz)(0)=1$)
$$
\frac{p_1\dots p_n\phi(0)}{K^{\alpha\mu^{-(n+1)}}(0)}\to 1,\quad n\to\infty.
\tag 6.24
$$
Therefore,
$$
\lim_{n\to\infty}
\frac{p_1\dots p_n\phi(0)}{K^{\alpha\mu^{-(n+1)}}(0)}
\frac{K^{\alpha\mu^{-n}}(0)}{p_1\dots p_{n-1}\phi(0)}=
\lim_{n\to\infty}\frac{ p_n}{\Cal P({\alpha\mu^{-(n+1)}})}=1.\tag 6.25
$$
It is a characteristic property of domains of Widom type, that
$$
0<\inf_{\alpha\in\Gamma^*} K^\alpha(0)\le
\sup_{\alpha\in\Gamma^*} K^\alpha(0)<\infty.
$$
Therefore, (6.25) implies
$$
p_n-\Cal P(\alpha\mu^{-(n+1)})\to 0,\quad n\to\infty.\tag 6.26
$$
Recall, that $K^\alpha(0) K^{\alpha^{-1}\mu}(0)=b'(0)$, hence
$$
\Cal P(\alpha)=\frac{K^\alpha(0)}{K^{\alpha\mu}(0)}=
\frac{K^{\alpha^{-1}}(0)}{K^{\alpha^{-1}\mu}(0)}=\Cal P(\alpha^{-1}).
$$
So we can rewrite (6.26) in the form (6.22).

To prove (6.23) we use (6.4). Since,
$$
\align
 &P_n b^n-\frac{\phi}{(- z')b^2}K^{\alpha^{-1}\mu^{n+2}}\to 0,\\
 &(zb)P_{n-1} b^{n-1}-\frac{\phi}{(-z')b^2}(zb)K^{\alpha^{-1}\mu^{n+1}}\to 0,\\
 &(p_n-\Cal P(\alpha^{-1}\mu^{n+1}))\frac{\phi}{(-
z')b^2}K^{\alpha^{-1}\mu^{n+2}}
 \to 0,\quad n\to\infty,
 \endalign
 $$
 we have
 $$
(zP_{n-1}-p_n P_n) b^n-\frac{\phi}{(- z')b^2}
(zb K^{\alpha^{-1}\mu^{n+1}}- \Cal
P(\alpha^{-1}\mu^{n+1})K^{\alpha^{-1}\mu^{n+2}})\to 0,\quad n\to\infty.
 $$
 Then, due to the maximum principle for analytic functions,
 $$
(zP_{n-1}-p_n P_n) b^{n-1}-\frac{\phi}{(-z')b^3}
(zb K^{\alpha^{-1}\mu^{n+1}}- \Cal
P(\alpha^{-1}\mu^{n+1})K^{\alpha^{-1}\mu^{n+2}})\to 0,\quad n\to\infty
 $$
 (note, that this function has no pole at the origin). Using the recurrence
relations and putting
 $\zeta=0$, we get
 $$
 \left\{\frac{q_{n-1}}{p_1\dots p_{n-1}}-\frac{\phi(0)}{b'(0)}\Cal
 Q(\alpha^{-1}\mu^{n+1})K^{\alpha^{-1}\mu^{n+1}}(0)\right\}\to 0,\quad
n\to\infty.
 $$
 Using again the identity $K^\alpha(0) K^{\alpha^{-1}\mu}(0)=b'(0)$ and
 (6.24), we get (6.23).\qed
\enddemo


\Refs
\ref\no 1
\by N.I. Akhiezer
\paper Orthogonal polynomials on several intervals.
\jour Soviet
Math. Dokl. 
\vol 1 
\yr 1960 
\pages 989--992
\endref


\ref\no 2
\by {N.I. Akhiezer}
\book {Elements of the theory of elliptic functions}
\publ {American Mathematical Society}, {Providence, RI}
\yr {1990}
\endref

\ref\no 3
\by N.I. Akhiezer and Yu.Ya. Tomchuk
\paper
On the theory of orthogonal polynomials over several intervals
\lang Russian
\jour Dokl. Akad. Nauk SSSR 
\vol 138 
\yr 1961
\pages 743--746
\endref

\ref \no 4
\by {A.I. Aptekarev}
\paper {Asymptotic properties of polynomials orthogonal on a system of
             contours, and periodic motions of {T}oda chains}
\jour {Mat. Sb. (N.S.)}
\vol {125(167)}
 \yr {1984}
   \pages  {231--258}
   \transl English transl. in  Math. USSR Sb., 53 (1986), 233--260
   \endref
   
   \ref\no 5
   \by S.N. Bernstein
   \paper Sur les polynomes orthogonaux relatifs \'a un segment fini, I
   \jour Journ. de math. pures et appl.
   \vol 9\pages 127--177\yr 1930
   \endref

\ref\no 6
\by L. Carleson
\paper On $H^\infty$ in multiply connected domains
\inbook Conference on harmonic analysis in honor Antoni
Zygmund. (eds. W Beckner, {\it et al.}). vol. II.
\publ Wadsworth
\yr 1983
\pages 349--372
\endref

\ref\no 7
\by R. Carmona and J. Lacroix
\book Spectral Theory of Random Schr\"odinger Operators 
\publ Birkh\"auser
\yr 1990
\endref

\ref\no 8
\by W. Craig 
\paper The trace formula for Schr\"odinger operators on the line
\jour Commun. Math. Phys.
\vol 126
\pages 379--408
\yr 1989
\endref

\ref\no 9
\by H.L. Cycon, R.G. Froese, W. Kirsch, B. Simon
\book Schr\"odinger Operators with Application to Quantum Mechanics and
Global Geometry
\publ Springer--Verlag, Berlin
\yr 1987
\endref

 

\ref\no 10
\by P. Deift
\book Orthogonal polynomials and random matrices:
a Riemann--Hilbert approach
\publ Courant Lecture Notes in Math., Vol. 3,
Courant Institute of Math. Sciences
\yr 1999
\endref

\ref\no 11
\by {B.A. Dubrovin, I.M. Krichever and S.P. Novikov}
 \inbook {Dynamical systems. {I}{V}}
            \publ {Springer-Verlag}, {Berlin}
\yr {1990}\pages  173--280
 \endref



\ref\no 12
\by J. Garnett
\book Bounded analytic functions
\publ Academic press
\yr 1981
\endref


\ref\no 13
\by J.S. Geronimo and W. Van Assche
\paper Approximating the weight function for orthogonal polynomials on
several intervals
\jour J. Approx. Theory
\vol 65
\yr 1991
\pages 341--371
\endref

\ref\no 14
\by M. Hasumi
\book Hardy Classes on Infinitely Connected Riemann Surfaces
\publ Lecture Notes in Math. 1027, Spinger Verlag, Berlin and New York 
\yr 1983
\endref

\ref\no 15
\by P. Jones
\paper Some problems in complex analysis
\inbook The Bieberbach Conjecture. Proc. of the Symp. on the Occasion
of the Proof (eds., A. Baernstein and D. Drasin)
\publ Amer. Math. Soc., Providence, RI
\yr 1986
\pages 105--108
\endref

\ref\no 16
\by P. Jones and D. Marshall
\paper
Critical points of Green's function, harmonic measure,
and the Corona problem
\jour
Arkiv f\"or Matematik
\vol 23\pages 281--314\yr 1985
\endref

 
\ref\no 17
\by D.S. Lubinsky
\book Asymptotics of Orthogonal Polynomials: some old,
some new, some identities
\publ to appear in Proc. of the Conf. on
"Rational Approximation", (eds., A. Cuyt, B. Verdonk),
Antwerp June 6--12, 1999, Kluwer Academic press
\endref

\ref\no 18
\by {A. Magnus}
 \paper {Recurrence coefficients for orthogonal polynomials on connected
             and nonconnected sets}
\inbook {Pad\'e approximation and its applications (Proc. Conf., Univ.
             Antwerp, Antwerp, 1979)}
    \pages {150--171}
\publ {Springer}, {Berlin}
 \yr {1979}
 \endref
 
\ref\no 19
\by V.A. Marchenko
\book Sturm--Liouville Operators and Applications 
\publ Birkh\"auser Verlag
\yr 1986
\endref


\ref\no 20
\by R. Nevanlinna
\book Analytic Functions
\publ  Springer Verlag, Berlin  
\yr 1970
\endref

\ref \no 21
\by E.M. Nikishin
\paper The discrete Sturm--Liouville operator and
some problems of function theory
\lang Russian
\jour Trudy Sem. Petrovsk. 
\vol 10 
\yr 1984
\pages 3--77
 \transl English transl. in Soviet Math., 35 (1987), 2679--2744
\endref

\ref\no 22
\by {E.M. Nikishin and V.N. Sorokin}
\book {Rational approximations and orthogonality}
 \publ {American Mathematical Society}, {Providence, RI}
 \yr {1991}
 \endref

 

\ref\no 23
\by L. Pastur and A. Figotin
\book  Spectra of Random and Almost--Periodic Operators
\publ  Springer Verlag, Berlin  
\yr 1986
\endref


\ref\no 24
\by {F. Peherstorfer}
 \paper {Elliptic orthogonal and extremal polynomials}
 \jour {Proc. London Math. Soc. (3)}
  \vol {70}
  \yr {1995}
  \pages {605--624}
   \endref
   
   \ref\no 25
 \by {F. Peherstorfer}
 \paper On Bernstein--Szeg\"o orthogonal polynomials on several intervals
 \jour SIAM J. Math. Anal.
 \vol 21
 \yr 1990
 \pages 461--482
 \endref
 
 
   \ref\no 26
 \by {F. Peherstorfer}
 \paper On Bernstein--Szeg\"o orthogonal polynomials on several intervals II:
 orthogonal polynomials with periodic reccurence coefficients
 \jour J. Approx. Theory 
 \vol 64
 \yr 1991
 \pages 123--161
 \endref
 
 
   \ref\no 27
 \by {F. Peherstorfer and R. Steinbauer}
 \paper {On polynomials orthogonal on several intervals}
  \jour {Ann. Numer. Math.}
   \vol {2}
     \yr {1995}
   \pages {353--370}
  \endref
  
  
\ref\no 28
\by Ch. Pommerenke
\paper
On the Green's function of Fuchsian groups
\jour
Ann. Acad. Sci. Fenn.
\vol 2\pages 409--427\yr 1976
\endref


\ref\no 29
\by M. Sodin and P. Yuditskii
\paper
Almost periodic Jacobi matrices with homogeneous spectrum,
infinite dimensional Jacobi inversion,
and Hardy spaces of character--automorphic functions
\jour Journ. of Geom. Analysis 
\vol 7\pages 387--435\yr 1997
\endref

\ref\no 30
\by G. Szeg\"o
\book Orthogonal polynomials
\publ 4th ed., Amer. Math. Soc. Colloq. Publ.,
Vol. 23, Amer. Math. Soc., Providence, R.I.
\yr 1975
\endref

  \ref\no 31
     \by { M. Toda}
    \book {Theory of nonlinear lattices}
   \publ {Springer-Verlag}, {Berlin}
    \yr {1989}
    \endref
   
\ref\no 32
\by Yu. Tomchuk 
\paper
Orthogonal polynomials over a system of intervals on the real line 
\lang Russian
\jour Zap. Fiz.-Mat. Fak. i Khar'kov. Mat. Ob\v s\v c.
(4) 
\vol 29 
\yr 1963
\pages 93--128
\endref

\ref \no 33
\by W. Van Assche
\book Asymptotics for Orthogonal Polynomials
\publ  Lecture Notes in Math. 1265,
 Springer, Berlin
\yr 1987
\endref

\ref\no 34
\by H. Widom
\paper Extremal polynomials associated with a system of curves in the
complex plane
\jour Adv. Math.
\vol 3
\yr 1969
\pages 127--232
\endref

\ref\no 35
\by H. Widom
\paper
The maximum principle for multiple valued analytic functions
\jour
Acta Math.
\vol 126
\pages 63--81
\yr 1971
\endref

\ref\no 36
\by P. Yuditskii
\paper Two remarks on Fuchsian groups of Widom type
\publ Oper. Theory Adv.  Appl.\toappear
\endref

\endRefs
\enddocument